\documentclass[journal]{IEEEtran}
\usepackage[english]{babel}
\usepackage[utf8]{inputenc}
\usepackage[ruled, vlined,linesnumbered]{algorithm2e}
\usepackage{authblk}
\usepackage[numbers,sort&compress]{natbib}
\usepackage{pdflscape}
\usepackage[space]{grffile}
\usepackage[cmex10]{amsmath}
\usepackage{amsthm}
\usepackage{amssymb}
\usepackage{xspace}
\theoremstyle{definition}
\newtheorem{definition}{Definition}[section]

\usepackage{booktabs}
\usepackage[inline]{enumitem}
\usepackage{tabularx}
\renewcommand{\vec}[1]{\ensuremath{\mathbf{#1}}}

\DeclareMathOperator*{\argmin}{arg\,min}
\DeclareMathOperator{\var}{Var}
\newcommand{\abs}[1]{\ensuremath{\lvert#1\rvert}}

\providecommand{\orcidID}[1]{}

\usepackage{etoolbox} 
\BeforeBeginEnvironment{appendices}{\clearpage}
\usepackage{url}
\usepackage[allcolors=blue,bookmarks,colorlinks]{hyperref}
\newcommand{\khd}{$k$-\emph{HD}\xspace}
\newcommand{\thd}{$\tau$-\emph{HD}\xspace}
\newcommand{\khdr}{$k$-\emph{HDR}\xspace}
\newcommand{\thdr}{$\tau$-\emph{HDR}\xspace}

\newcommand{\assign}{\ensuremath{\gets}}
\newcommand{\mDM}{\ensuremath{m_\textsc{dm}}}
\newcommand{\Fdm}{\ensuremath{F_\textsc{dm}}}

\begin{document}
\bstctlcite{BSTcontrol}
\title{Detection of Hidden Objectives and Interactive Objective Reduction}
\author[1]{Seyed Mahdi Shavarani}
\author[1]{Manuel López-Ibáñez}
\author[1]{\IEEEmembership{Member,~IEEE}~ Richard Allmendinger}
\affil[1]{Alliance~Manchester~Business~School,~University~of~Manchester, Manchester,~M13~9PL,~UK \authorcr Email: {\tt \{seyedmahdi.shavarani,manuel.lopez-ibanez, richard.allmendinger\}@manchester.ac.uk}\vspace{1.5ex}
\authorcr This work has been submitted to the IEEE for possible publication. Copyright may be transferred without notice, after which this version may no longer be accessible.
}

\maketitle

\begin{abstract}
In multi-objective optimization problems, there might exist \emph{hidden} objectives that are important to the decision maker but are not being optimized. On the other hand, there might also exist \emph{irrelevant} objectives that are being optimized but are of less importance to the DM. The question that arises here is whether it is possible to detect and reduce irrelevant objectives without deteriorating the quality of the final results? In fact, when dealing with multi-objective problems, each objective implies a significant cost best avoided if possible. However, existing methods that pertain to the reduction of objectives are computationally intensive and ignore the preferences of the decision maker. In this paper, we propose an approach to exploit the capabilities of interactive evolutionary multi-objective optimization algorithms (EMOAs) and the preference information provided by the decision maker, to detect and eliminate the irrelevant objectives during the optimization process and replace them with hidden ones, if any. The proposed method, which is based on uni-variate feature selection, is computationally effective and can be integrated into any ranking-based interactive EMOA. Using synthetic problems developed in this study, we motivate different scenarios in the experiments and prove the effectiveness of the proposed method in improving the computational cost, the total number of objective evaluations and the quality of the final solutions.

\end{abstract}

\begin{IEEEkeywords}
Interactive Multi-Objective Optimization, Hidden Objective Functions, Redundant Objectives, Irrelevant Objectives, Machine Learning, Dimension Reduction, Feature Selection
\end{IEEEkeywords}

\section{Introduction}\label{introduction}
Interactive Evolutionary Multi-Objective Algorithms (iEMOAs) mostly assume that the objectives being optimized match the ones that interest all potential Decision Makers (DMs), however, this may not be the case. As instance, the DM may not be consulted during the modelling phase or his desires may not be well satisfied during the group decision-making on final objectives to be optimized. As the result, there may exist objectives that are optimized but the DM is not really interested in and on the other hand, there are objectives that the DM cares for but are not being optimized.
These inconsistencies may cause the elicited preference information  appear non-rational, e.g. when a solution that is dominated with respect to the modelled objectives is preferred by the DM over a non-dominated one because the former is better than the latter with respect to features that are not modelled as objectives by the algorithm. 
A real-world example would be the case discussed by \citet{RamMirSer2021foodplan}, where a problem is defined to optimize the food served to school children in terms of not only cost, but also a number of metrics of nutritional value and food variety. The variations in food types is optimized by penalizing repetitions of each food group. Repetitions in each food group may be considered as an independent objective. However, in that study they are merged ``a priori'' into a single objective. A DM that is interacting with the system may consider a given non-optimal menu (not) appealing but may not be able to justify this decision. Is it because there is a hidden tendency to optimize a certain type of food type? When making decisions, does the DM really considers all nutritious and food variety metrics as assumed by the model? Does the system really need to optimize all of the objectives? Or are there any extra objectives that the DM considers but are not reflected in the model?  Another example would be optimizing the  behavior of robots~\cite{TriLop2015plos}. There might exist a set of potential objectives (e.g., average speed, morphology, accuracy, etc.) for how well the robot  behaves but, at the end, a human decides which behaviors are better by observing the robot perform tasks. Thus, there might be a disconnection between the set of objectives being optimized and the set of (relevant) objectives that the DM uses in evaluating alternative behaviors. 

Under such circumstances, it is desirable to have an algorithm that identifies relevant and irrelevant objectives and updates the optimization model in such a way that only the relevant objectives take part in the optimization process.
Removing irrelevant objectives is specifically important in many-objective optimization problems. In some scenarios there is a large number of solution features that DMs may wish to optimise and there is no way but to model the problem as a many-objective optimization problem. However, if not all objectives are relevant to a particular DM, removing those irrelevant objectives will be beneficial to solving many-objective problems where the non-dominance strategy loses its selection pressure as most of the solutions become non-dominated with respect to the current population \cite{KhaYaoDeb2003,IshTsuNoj08,BroZit2006allobjectives}.
 Even in the case of many-objective algorithms that overcome problems associated with the dominance-based EMOAs \cite{ZitKun2004ppsn,BeuNauEmm2007ejor}, objective reduction techniques are still beneficial because, first, the running time increases exponentially with the number of objectives \cite{BeuFonLopPaqVah09:tec} and, second, the number of non-dominated solutions required for covering the whole Pareto front increases exponentially with the number of objectives~\cite{PurFle2003cec, SinIsaTap2011pareto}, hence, complicating the decision-making phase \cite{BroZit2009ec}.
 
Given the benefits of reducing the number of objectives, previous research has considered removing objectives that are highly correlated to other objectives \cite{ BroZit2009ec, SinSaxDeb2013ASC}, or their elimination does not have much of an effect on dominance relations~\cite{brozit2006dimensionality, BroZit2006allobjectives}.Here, we consider an approach for objective reduction that is different and complementary to previous ones.

Here, for the first time in the literature, we consider the preference information elicited  in interactive EMOAs as an opportunity to detect irrelevant objectives during the optimization  and update the model in a way that only relevant objectives, i.e. important for the DM, are optimized, while the rest of the objectives are made inactive. The relevance of objectives may change from a DM to another. The proposed approach gives the interactive EMOA the ability to adapt to different DMs and update objectives that are relevant to the DM who is interacting with the EMOA. The main contributions of this study can be summarized as follows:
\begin{itemize}
    \item A motivation and formal definitions of the problem of identifying relevant objectives amongst a set of potential objectives in interactive multi-objective optimization.
    \item Proposal of an approach to detect relevant and hidden objectives and to update the set of objectives during the optimization to optimize only the relevant ones. The approach draws on feature selection methods and can be applied to any ranking-based interactive EMOA. 

    \item  An approach is proposed to simulate hidden and irrelevant objectives which can be applied to any multi-objective problem and utility function. This approach is demonstrated for DTLZ problems \cite{DebThiLau2005dtlz} and  multi-objective NK-landscapes $\rho$MNKs~\cite{VerLieJou2013ejor}.
    \item A validation of the proposed detection method is performed which includes: 
    \begin{enumerate}[label=\roman*.]
        \item Problems of varying dimensionality, complexity, and Pareto front structure.
        \item Different utility functions to represent different DMs.
        \item Different feature selection methods to detect relevant objectives, and
        \item A sensitivity analysis to understand the performance impact of key parameters of the proposed approach.
    \end{enumerate}
\end{itemize}

The results of the experiments show that our proposed method is able to quickly replace irrelevant objectives with relevant ones and significantly improve the quality of final solutions. Moreover, it is observed that, even in cases where the utility of final solution is not favored by the proposed method, the computational effort is highly reduced.
The rest of the paper is organized as follows. The problem is defined and formulated in Section \ref{definitions}. A summarized background on previous efforts towards objective reduction is given in Section~\ref{background}. In Section~\ref{methods} the proposed solution method are elaborated in detail. The experimental setup is laid out in Section~\ref{experiments}. The results of the experiments are discussed in Section~\ref{results}.
Finally Section~\ref{conclusion} provides conclusions and future research directions.

\section{Definitions and formulations} \label{definitions} Let us consider a
many-objective optimization problem in which we could potentially optimize
simultaneously a set $\{f_1,\dotsc,f_m\} \in F$ of $m$ objectives:
\begin{equation}\label{eq:mo}
  \begin{split}
    \text{Minimize}\quad&\bigl( f_1(\vec{x}), \dotsc, f_m(\vec{x})  \bigr)\\
    \text{subject to}\quad&\vec{x} \in \mathcal{X}\\
  \end{split}
\end{equation}
where each objective $f_i\colon \mathcal{X} \to \mathbb{R}$ depends on a solution vector $\vec{x}=(x_{1},\dotsc,x_n)$ of
$n$ design (or decision) variables, and $\mathcal{X}$ is the feasible decision space.

\begin{definition}[Domination]
A solution $\vec{x}\in \mathcal{X}$ dominates solution $\vec{y} \in \mathcal{X}$ if  $\vec{f}(\vec{x})$ is not worse than  $\vec{f}(\vec{y})$ in any objective and $\vec{f}(\vec{x})$ is strictly better than $\vec{f}(\vec{y})$ in at least  one objective.
\end{definition}

\begin{definition}[Non-domination]
A feasible solution $\vec{x} \in \mathcal{X}$ is said to be non-dominated if there is no solution $\vec{y}\in \mathcal{X}$ that dominates it. Non-dominated solutions are also known as Pareto optimal solutions.
\end{definition}

\begin{definition}[Pareto front]
The projection of the Pareto solutions on the objective space is known as the Pareto front (PF).
\end{definition}

\begin{definition}[Potential objectives]
All $m$ objectives in $F$ are called \emph{potential} objectives. 
\end{definition}

\begin{definition}[Redundant objectives] An objective is called \emph{redundant} if it can be eliminated without changing the set of Pareto optimal solutions \cite{GalLeb1977ejor, AGRELL1997ejor}.
\citet{SaxRayDeb2009constrained} extend this definition to include objectives that are not conflicting with a non-redundant objective.
  \end{definition}

The above definitions are independent of the preferences of a DM interacting with the EMOA. Interactive methods can be classified into ad-hoc and non-ad-hoc methods~\cite{Steuer1986}. Non-ad-hoc methods assume there exists a utility function (UF) which derive the DM's decisions but are unknown to the EMOA~\cite{SteGar1991computational}. Due to the popularity of UFs in modeling preference models, vast majority of iEMOAs are non-ad-hoc methods thus we focus on those in the remainder of the paper.

Without loss of generality, here we assume 
an unknown UF $U\colon \mathbb{R}^{\mDM} \to [0,1]$ the value of which depends mainly on
$\Fdm \subseteq F$ objectives ($\mDM=\abs{\Fdm} \leq m = \abs{F}$). An instance will be a linear additive utility function where the weights associated to irrelevant objectives are close to zero.

\begin{definition}[Irrelevant objectives]
An objective is ``\emph{irrelevant}'' to the DM if it does not
appear in the DM's utility function. In other words, irrelevant objectives ($F\setminus\Fdm$) do not contribute to the utility of solutions.
\end{definition}

Interactive methods try to somehow estimate the DM's UF by exploiting the preference information that is elicited from the DM in the interactions that are performed during the optimization process. The estimated UF is then used in the next generations of the interactive EMOA  to guide the search towards that part of the Pareto front that is of interest to the DM. Let us assume that, at some moment of its execution, a given EMOA only optimizes a subset $\hat{F} \subseteq F$ ($\hat{m}=\abs{\hat{F}} \leq m$) of the potential
objectives. 

\begin{definition}[Active objective]
  An objective is \emph{active} if it is optimized by the EMOA. The set of
  active objectives is represented by $\hat{F}$. Inactive objectives
  ($F\setminus\hat{F}$) are either evaluated but ignored by the EMOA or not
  evaluated at all, unless required by the DM during an interaction.
\end{definition}

\begin{definition}[Objective evaluations \& computational efficiency]
An objective evaluation is defined as a single evaluation of any of the objectives corresponding to a solution ($f_i(\vec{x})$) and is different than solution evaluation (\vec{f}(\vec{x})) which corresponds to evaluation of all of the objectives of a given solution. In line with our previous definitions, each solution evaluation is equivalent to $\hat{m}$ objective evaluations. On the other hand, computational cost can be defined in a broader scale inclusive of objective evaluations and other
aspects of the solution procedure such as the cost of non-domination sorting, which are again affected by the number of active objectives which participate in sorting and optimization. Thus, we use the number of active objectives and objective evaluations as an indicator of computational efficiency; i.e. reducing the number of active objectives reduces computational costs.
\end{definition}

\begin{definition}[Hidden objectives]
  An objective $f_i \in F$, is said to be \emph{hidden} if $f_i \in \Fdm \land f_i \notin \hat{F}$. That is, an objective that appears in the DM's UF (relevant), but is not (currently) active.  If $\Fdm \subset \hat{F}$, then no hidden objectives exist, but the EMOA is still optimizing some irrelevant objectives. Similarly,  if $\Fdm = \hat{F}$, then neither hidden nor irrelevant objectives exist, and the EMOA is optimizing precisely the objectives that the DM cares about.
\end{definition}

Being able to reveal hidden objectives is beneficial, since otherwise there is
a disconnection between what the algorithm is optimizing and what the DM cares
about. Moreover, the interactive algorithm may become confused if the
interaction with the DM is consistent with the dominance criterion  for the
objectives in $\Fdm$ but appears to contradict it for the objectives in
$\hat{F}$.

Without any prior information about $\Fdm$, most EMOAs would optimize all
potential objectives \mbox{$\hat{F} = F$} leading to a many-objective
optimization problem, which is inherently challenging for EMOAs as discussed
above.

From these definitions, it can be concluded that while \emph{redundant} objectives are determined based on the structure of the problem, \emph{irrelevant} and \emph{hidden} objectives are defined from the DM's perspective.  While there are studies on detection and elimination of redundant objectives, which we review in the next section, there is no prior research on the identification of irrelevant and hidden objectives to the best of our knowledge.  Our focus here is to fill this gap and we propose a method to tackle it in Section \ref{methods}.

\section{Background}\label{background}
 A brief review of the literature on the reduction of \emph{redundant} objectives is presented here. 
 These methods  are applied either ``a priori'' to facilitate the optimization process or ``a posteriori'' to assist the DM in decision making.
Presented studies either consider correlations among objectives or their importance in preserving the dominance relations among solutions.  We have categorized  existing methods accordingly in such a way that makes it easy to compare similar studies and also discuss their limitations. We acknowledge there are reduction methods that attend to the dimensions of the decision space, such as those discussed in \cite{AllKno2010variables, ZHAZhaChen2020surrogate, ChuZhaFu2015HESS}, as well as those that focus on constraint reduction \cite{SaxDeb2008Dimensionality,SaxDeb2007Trading}. However, they fall out of the scope of this research as our focus is on the objective space. Many of these studies use the term \emph{dimension} reduction for the same concept. However, to avoid confusion with methods that reduce the number of decision variables~\cite{AllKno2010variables}, we use the term ``objective reduction''.

\subsubsection{Subset selection \& preserving Pareto front solutions}
\citet{GalLeb1977ejor} proposed one of the earliest  approaches that considered preserving the PF. Although the method was later extended by \citet{AGRELL1997ejor}, they both make exhaustive assumptions about the problem structure that are impossible to validate for real-life problems.
\citet{brozit2006dimensionality, BroZit2006allobjectives} 
consider objective reduction while preserving the PF. They introduce two minimum
objective subset selection (MOSS) methods: $\delta$MOSS finds a minimal set of objectives with maximum error of $\delta$ and $k$-EMOSS finds a subset with fixed size $k$ having the minimum error. Error here is a measure of change in the set of PF solutions. The performance of these methods is compared with PCA-based methods in \cite{BroSaxDeb2008handling} and the results indicated the competitiveness of the methods. In \cite{BroZit2007hypervolumeReduction}, the authors use their minimal set approach online in simple indicator-based evolutionary algorithm (SIBEA) \citep{ZitBroThi2007emo}. As MOSS problems are computationally expensive, \citet{Gui2011objred} proposed a linear programming model for $k$-EMOSS problem, which was demonstrated to have acceptable performance.  

\citet{JaiCoeCha2008} use an unsupervised feature selection technique~\cite{MitMurPal2002unsupervised} for preserving the most conflicting objectives and identifying the subset of non-redundant objectives. In their approach, neighborhoods of objective functions are formed based on correlation of some non-dominated solutions. Then an objective is selected from each compact neighborhood while the rest are dropped.
Later, this technique was incorporated into NSGA-II~\cite{Deb02nsga2} to construct an online reduction method~\cite{JaiCoeUri2009onlinered}.

\citet{SinIsaTap2011pareto} focus on the Pareto corner points to identify which objectives are sufficient to reproduce the PF. Some drawbacks of their approach are: not capturing the whole PF, losing objectives that are contributing somewhere on the PF away from corner points and possible overestimation of dimensions.

\subsubsection{Recombining non-conflicting objectives}    

    \citet{FreFleGui2013NonParametric} used the harmonic level, introduced in~\cite{PurFle2007tec, PurFle2003ConflictHarmony, brozit2006dimensionality, giaPurFle2013decompositon}, to identify objectives that can be merged into a new compound scalarized objective with minimal effect on the PF. As opposed to conflict, there is ``harmony'' among objectives when improvement of either of them does not lead to deterioration of the other \cite{PurFle2003ConflictHarmony, PurFle2007tec, giaPurFle2013decompositon}. When there is harmony between two objectives, they can be merged without any loss in quality of the results. Similar to other objective reduction methods, the DM preferences are not considered in the reduction process in \cite{FreFleGui2013NonParametric}. However, the DM can decide on the final number of objectives ``a priori''.  
    Similarly, \citet{FreFleGui2015aggregation} use aggregation trees to identify harmonious objectives that could be aggregated into single scalarized objectives with minimal loss.
    
\subsubsection{PCA \& Correlation among objectives}
One of the first methods that relied on linear PCA to spot redundant objectives was by \citet{debsax2006searching} that was implemented within the NSGA-II framework. Linear PCA relies on  correlation among objectives to detect conflicting objectives, thus it does not guarantee that the dominance relation is preserved \cite{BroZit2009ec}.Furthermore, linear PCA does not provide a good performance for non-linear data, is limited by the type of the data it can handle and it might fail identify all non-conflicting objectives properly, but in local regions where solutions are close it may still be used for determination of conflicts and harmony between objectives \cite{LygMarFle2010Reduction}. To address these issues, two more  methods were introduced in \cite{saxdeb2007non-linear}, which used correntropy PCA \cite{XuPokPai2006correntropyPCA} and maximum variance unfolding as was previously used in \cite{WeiShaSau2004dimension, weiSau2006introduction}. PCA was also used in \cite{GoeVaiHaf2007response}, however the method was validated on a problem with $4$ objectives only. 
PCA was also used inline with NSGA-II in \cite{PozRui2012PCAreduction, SaxDurDebZha2013pca} among others,  to handle redundant objectives in non-linear contexts.
PCA-based methods are applied on non-dominated solutions and thus their performance relies heavily on proper approximation of the PF \cite{SinIsaTap2011pareto}. 

In their study, \citet{CosOli2010biplotsReduction} have demonstrated that objectives that are deemed redundant by PCA may be ``informative'', i.e.,  contain trade-off information that would be lost if omitted. They suggest using PCA coupled with BIPLOT~\cite{Gabriel1971biplot} representation to confirm if a PCA-identified redundant objective is not ``informative''. The applicability of the method is limited to the problems with only one independent objective.
  
\subsubsection{Visualization}
Some studies suggest identification of redundant objectives through visualization \cite{ObaSas2003visualization, MatMes2003AIAA, KopYos2007visualization, FieEve2013visualising}. However, most of these studies aim to facilitate intuitive visualization of the solutions rather than the optimization process itself and are limited in the number of objectives they can handle. Thus, they use objective reduction methods such as PCA to enable the visualization and hence bring new insights to the DM.
\\

\section{Methods} \label{methods}
As laid out, existing approaches for objective reduction suffer from dependency on non-dominated solutions, expensive computations, elimination of informative objectives, offline reduction, and not exploiting the DM's preferences. In this section, we illustrate a method that is able to use the DM's preferences, elicited during optimization by an interactive method, to identify irrelevant objectives as well as hidden ones dynamically.The experimental study will show that the proposed method is able to efficiently reduce the number of objectives being optimized by the optimizer, hence dedicating more computational resources to the optimization of relevant objectives only.

In a nutshell, our proposal works as follows. 
At some point during the run of the EMOA, the method interacts with the DM by showing the value of all potential objectives of a selected subset of solutions and asking the DM to rank the solutions. Feature selection applied to the rankings and the objective values identifies which objectives have the most significant effect on the ranking. The method uses this information to possibly activate currently inactive objectives and/or deactivate currently active ones. The EMOA then continues the search using the new set of active objectives. In what follows the method is described in detail.

\subsection{Feature Selection}\label{sec:feature_selection}
We explore two feature selection methods in this study: Uni-variate feature selection and Recursive Feature Elimination (RFE). This gives us the opportunity to investigate and compare the results with respect to different feature selection methods.
 Hereafter, feature and objective are used interchangeably in this context. 
\subsubsection{Uni-variate Feature Selection}
We propose the application of F-test uni-variate feature selection for identifying the most relevant features. In uni-variate methods, each feature is considered independently and any correlation between features is ignored \cite{DesMosMol2009FeatureSelection}.
 The uni-variate detection method applied here estimates the degree of linear dependency between two random variables. Let $T$ be the set of solutions presented to the DM at an interaction, where $\vec{f}_j \in T$ is the $j^\text{th}$ vector of objective values, and $f_{ji}$ denotes the value of its $i^\text{th}$ objective out of the $m$ potential objectives. The DM ranks the solutions according to her own preference of the objective values. The vector of rankings is given by $\vec{r}$, where $r_j$ is the rank corresponding to $\vec{f}_j \in T$.
 There is no restriction on the rankings and two solutions may have the same rank. 

 The procedure for F-test uni-variate feature selection can be described as follows (for reference on F-test feature selection see \cite{Heiman2000ResearchMethods}):

\begin{description}
\item[Step 1:] The correlation $\rho_i$ between each objective (feature) $i$ and $\vec{r}$ is computed as:
\begin{equation}\label{eq:rho}
  \rho_i=\sum_{j=1}^{|T|} \frac{(f_{ji}-\bar{f_{\cdot i}})\cdot(r_{j}-\bar{\vec{r}})}{\var(f_{\cdot i})  \var(\vec{r})}
  \end{equation}
  where $\bar{f_{\cdot i}}$ and $\var(f_{\cdot i})$ are, respectively, the mean and variance  of the $i^\text{th}$ objective over all solutions in  $T$;   and $\bar{\vec{r}}$ and $\var(\vec{r})$ are the same for the vector of rankings.

\item[Step 2:] The F-statistic for each objective is computed as:
    \begin{equation}\label{eq:f_statistic}
        \mathcal{F}_i=\dfrac{\rho_i}{1-\rho_i}\cdot(\abs{T}-2)
    \end{equation}
    Here, the F-statistic is a notion of how well an objective can explain the rankings provided by the DM.
    
  \item[Step 3:] The $p$-values corresponding to each F-statistic is calculated by any statistical software.
  \item[Step 4:] Features with lower $p$-value are selected. Number of selected features can be either fixed ($k=\hat{m}$) or variable. In the latter case, objectives with $p$-values less than a predetermined threshold $\tau$ are selected. These two variants are explained in Section~\ref{fixed vs variable}.
  \end{description}

The lower the $p$-value, the better is the corresponding objective function in explaining the DM's rankings. The pseudo-code of uni-variate feature selection is illustrated in Algorithm~\ref{alg:univariateFeatureSelection}.

\begin{algorithm}[tp]
  \caption{Uni-Variate Feature Selection}
  \label{alg:univariateFeatureSelection}
  \KwIn{
    \begin{description}
    \item[]$F$ : Set of all potential objectives
    \item[]$T$ : Set of ranked objective vectors
    \item[]$\vec{r}$ : vector of ranks 
    \item[]Either 
    \begin{itemize}
      \item $k < \abs{F}$ (for fixed number of objectives) or
      \item $\tau$ (for variable number of objectives)
  \end{itemize} 
    \end{description}
  }
\KwOut{$\hat{F}$: Selected objectives}
\For{$i \assign 1$ \KwTo $\abs{F}$}{
     Step 1: Calculate $\rho_i$ using Eq.~\eqref{eq:rho}\\
     Step 2: Calculate $\mathcal{F}_i$ using Eq.~\eqref{eq:f_statistic}\\
     Step 3: Calculate $p_i$ ($p$-value) from $\mathcal{F}_i$
   }
   \eIf{Fixed number of objectives}{
   $\hat{F}$ \assign $k$ objectives with lowest $p$-value}{
   $\hat{F} \assign  \{ f_i \in F \mid p_i < \tau \}$\label{inalgo:selection}} 
 \Return{$\hat{F}$}
\end{algorithm}
\subsubsection{Recursive Feature Elimination}\label{Recursive fs}
RFE is different than uni-variate feature selection in that, here, first logistic regression is used to build a model based on all the features to predict the rankings. In the next step, the feature with the minimum contribution to the constructed model is excluded from the selected subset~\cite{GuyWesBar2002rfe}.

The two steps are repeated until a preset number of objectives $k$ is selected (fixed number of objectives) or the minimum contribution of the remaining objectives is greater than a threshold $\tau$ (variable number of objectives). There are several ways to measure contribution to the model ($\phi$) such as the coefficient of the objective in the model or the importance of the objective. Here, we use the latter to be consistent with uni-variate feature selection. The importance of each feature is calculated based on the drop in the accuracy of the model when that feature is eliminated. 
The algorithm for RFE is depicted in Algorithm \ref{alg:RecursiveFeatureSelection}.

\begin{algorithm}[tp]
  \caption{Recursive Feature Elimination}
  \label{alg:RecursiveFeatureSelection}
  \KwIn{
  \begin{description}
  \item[]$F$: Set of all potential objectives
  \item[]$T$: Set of ranked solutions
  \item[]$\vec{r}$: vector of ranks
  \item[]Either 
  \begin{itemize}
      \item $k < \abs{F}$ (for fixed number of objectives) or
      \item $\tau$ (for variable number of objectives)
  \end{itemize} 
  \end{description}\\\\ 
  }
\KwOut{$\hat{F}$: Selected objectives}
$\hat{F} \assign F$\\
\While{True}{
  Step 1: $M$ \assign Build\_Model($T$, $\vec{r}$, $\hat{F}$)\\
 Step 2: $f_j \assign \argmin_{f_i \in \hat{F}} \phi(f_i)$\\ 
 \If{$\phi(f_j)>\tau$ or $|\hat{F}|=k$ or $|\hat{F}|=2$}{\textbf{break}}
 Step 3: $\hat{F} \assign \hat{F}\setminus f_j$
 }
 \Return{$\hat{F}$}
\end{algorithm}
\subsection{Fixed versus Variable Number of Active Objectives}\label{fixed vs variable}
The number of features (active objectives) selected can be defined in different ways. Here we explore the following two alternatives: 
\subsubsection{Fixed number of objectives ($k$)}
The optimization starts with
$k$ active objectives and this number is kept constant throughout the
optimization process such that activating an inactive objective implies
deactivating an active one. 

\subsubsection{Variable number of objectives}
We select the subset of objectives which meet a
predetermined threshold $\tau$.  In case of uni-variate feature-selection, the lower the value of $\tau$, the lower would
be the number of objectives with acceptable $p$-values. If there is only one
objective with a $p$-value lower than $\tau$, the two objectives with lowest
$p$-values are selected instead.
For RFE, in each step an objective is selected only when its contribution to the model is greater than $\tau$. Thus, higher values of $\tau$ corresponds to less number of active objectives.
Having two feature selection methods and two approaches with fixed and variable number of objectives as explained above, we have four total variations defined as follows:
\begin{enumerate}
    \item \khd: Uni-variate feature selection with fixed number of objectives
    \item \thd: Uni-variate feature selection with variable number of objectives
    \item \khdr: RFE with fixed number of objectives
    \item \thdr: RFE with variable number of objectives
\end{enumerate}
The proposed methods can be applied to any ranking-based interactive EMOA for objective reduction and/or detection of hidden objectives in order to find the objectives that are relevant to the DM. 
Here, we will focus on extending BCEMOA~\cite{BatPas2010tec} with our proposed method to show how the method can be integrated with any ranking based algorithm. In what follows, the modified BCEMOA here called BCEMOA-HD, is explained in detail. 

\newcommand{\Ninteractions}{\ensuremath{N^\textup{int}}}
\newcommand{\Nexa}{\ensuremath{N^\text{exa}}}
\newcommand{\Usvm}{\ensuremath{U_\textsc{svm}}}

\subsection{BCEMOA-HD}
As discussed, our proposed HD/HDR method can be applied to any ranking based interactive method. In this study, BCEMOA is selected and equipped with detection of hidden objectives.
BCEMOA \cite{BatPas2010tec} is an interactive EMOA based on NSGA-II.
The algorithm starts with a population of randomly generated solutions (\emph{pop}). The population is evolved with NSGA-II for $\text{gen}_1$ generations. Next, at
each interaction step, the best \Nexa solutions are selected from the evolved population and ranked by the DM. The data along their ranks are then used to train a Support Vector Machine (SVM) model to learn a utility function (\Usvm). The evaluations of the learned UF replaces the crowding distance in the next generations. Further interactions with the DM provide additional
samples to re-train the SVM model and improve the predictions of the learned utility.

Similar to the original BCEOMA, the BCEMOA-HD algorithm, proposed here, starts with a set of active objectives $\hat{F}$. All inactive objectives not in $\hat{F}$ do not need to be evaluated during the optimization and do not participate in dominance ranking and evolution of the population. However, during the interactions, all objectives in $F$ are evaluated for the solutions that are presented to the DM. Immediately after each interaction, the feature selection method described in Section \ref{sec:feature_selection} is applied to the objective vectors and their rankings to identify relevant objectives and update $\hat{F}$. Consequently, the population should also be updated to be evaluated for $f_i \in \hat{F}$. SVM is also used to learn \Usvm  based on active objectives in updated $\hat{F}$  and their rankings. An overview of  BCEMOA-HD is shown in Algorithm~\ref{alg:bcemoa}.

As described above, compared to the original BCEMOA we have modified the algorithm in lines \ref{featureSelectionLine} and \ref{updatePop}, where feature selection is deployed and the set of active objectives is updated. 
Another modification was applied to the original BCEMOA in the selection of the best solutions presented to the DM.
When there is no variance in the values of some objective, for example, because its values are near-optimal, then their correlation with the rankings provided by the DM is almost zero, which would result in the feature selection considering that objective as irrelevant and replacing it with a different one. 
To address this issue and to preserve the elitism, the solution that was ranked best by the DM in the last interaction is always included in the next set of solutions presented to the DM by BCEMOA-HD. 

\begin{algorithm}[tp]
  \caption{BCEMOA-HD}
  \label{alg:bcemoa}
   \KwIn{\\$\Ninteractions\colon$ Total number of interactions\\ $\Nexa\colon$ Number of training examples per interaction\\ $\textit{pop}\colon$ population of solutions\\ $\textit{gen}_1$: Generations before first interaction\\ $\textit{gen}_i$: Generations between two interactions\\
   $F$: Set of potential objectives\\
   $\hat{F}$: Set of active objectives\\
   Either 
  \begin{itemize}
      \item $k < \abs{F}$ (for fixed number of objectives) or
      \item $\tau$ (for variable number of objectives)
  \end{itemize} 
   
  }
\KwOut{The most preferred solution}

 $\textit{T} \assign \{\}$ \\
 $\textit{r} \assign \{\}$ \\
 \textit{pop} \assign NSGA-II for $\textit{gen}_1$ generations (\textit{pop})\\
 \For{$1$ \KwTo $\Ninteractions$}{
   $\textit{T}_i \assign \texttt{select}\,\Nexa\,\texttt{solutions}$\\
   Evaluate solutions in $\textit{T}_i$ for all objectives in $F$\\
  $\textit{r}_i$ \assign \texttt{DM\_ranks}($\textit{T}_i$)\\
  $\textit{T} \assign \textit{T}_i$ $\cup$ \textit{T}\\
  \textit{r} \assign $\textit{r}_i$ $\cup$ \textit{r}\\
  $\hat{F}$ \assign \textit{feature\_selection(\textit{T}, \textit{r}, $F$, $\tau$ / $k$)}\\\label{featureSelectionLine}
  Evaluate \textit{pop} for $f_i \in \hat{F}$\\\label{updatePop}
  \Usvm \assign \texttt{train\_SVM}(\textit{$T$ with objectives $\in\hat{F}$}, $\vec{r}$)\\
  \texttt{Crowding\_Distance} \assign \Usvm \\
  \textit{pop} \assign {NSGA-II for $\textit{gen}_i$ generations, only considering $\hat{F}$}
 }
 \Return{\textit{Best (first) solution in the \textit{pop} ranked by non-domination sorting and \Usvm}}
\end{algorithm}

\section{Experimental Setup} \label{experiments}\label{exp_setup}
To evaluate the effectiveness of our proposed method, we design a set of experiments that can comprehensively cover different aspects of the problem of identifying hidden and irrelevant  objectives. To do so, we propose a possible way of simulating hidden and irrelevant objectives. Two sets of well-known benchmark problems are selected from the literature. Each problem features various difficulties to different parts of the solution method including convergence towards most preferred point and performance in highly rugged landscapes. In what follows, a detailed description of the design of the experiments is laid out.

In the experiments with variable number of objectives we seek to investigate how the method performs for objective reduction purposes. Thus, for these set of experiments all objectives are active from the start of the run ($\hat{F}=F$), while for fixed number of objectives only designated objectives are active ($\hat{F}\subset F$) at the beginning of the optimization.

\subsection{Synthetic problems with hidden objectives}
We create synthetic problems that feature irrelevant and hidden objectives by extending existing benchmark  problems as follows.
Given a problem with $m= \abs{F}$ potential
objectives, we extend it with a
binary vector $\vec{d} \in \{0,1\}^m$ that specifies which  objectives are considered by the optimizer, i.e.,
$d_i = 1$ iff $f_i \in \hat{F} \subset F$, where $f_i$ indicates the $i^\text{th}$ objective function. That is, given a solution $\vec{x}$, whose objective vector
is $\vec{f}(\vec{x}) = (f_1(\vec{x}), \dotsc, f_m(\vec{x}))$, the optimizer only considers $\vec{\hat{f}}(\vec{x})=\vec{f}(\vec{x}) \odot \vec{d}$, where $\odot$ denotes the
element-wise product of two vectors. The optimizer is able to change the set of active objectives by changing the  vector $\vec{d}$.

On the other hand, feature selection methods and the DM  have access to
$\vec{f}(\vec{x})$, that is,
the unknown UF used by the DM when asked to compare pairs of objective vectors  actually evaluates
$U(\vec{f}(\vec{x}))$.

\subsection{Benchmark Problems}

We consider two well-known numerical and binary benchmark problems, namely,  multi-objective NK landscape problems with correlation between objectives ($\rho$MNK) \cite{VerLieJou2013ejor} and DTLZ problems 
\cite{DebThiLau2005dtlz} 
 with $m \in \{4, 10, 20\}$ objectives. 
 Problems with $m=4$ help us to better understand and investigate the dynamics of the proposed methods, while larger number of objectives allows us to evaluate the efficiency of the feature selection with variable number of objectives in many-objective problems.

$\rho$MNK problems are used to analyse the effects of correlation among objectives and smoothness of the landscape on the performance of the proposed method. We consider $\rho$MNK instances  with different values of correlation among objectives $\rho \in \{-0.25, 0, 0.25, 0.5, 0.75, 0.9\}$, taking into account the restriction that  $\rho \geq -1/(m-1)$~\cite{VerLieJou2013ejor} and different values of parameter $K$, which controls the smoothness of the landscape, namely, $K \in \{1, 4, 6, 8\}$ for problems with 4 objectives and $K \in \{1, 5, 10, 15\}$ for many objective problems, considering the constraint $K<n$. The greater the value of $K$, the more rugged is the fitness landscape. The value of $n$ is kept fixed at 10 for problems with $m=4$, 20 for problems with $m=10$ and 30 for problems with $m=20$ for $\rho$MNK problems.

We use DTLZ1, DTLZ2, DTLZ7 from DTLZ test suit, which were also used in the experiments on  BCEMOA by its authors \cite{BatPas2010tec} and also in \cite{BroZit2007hypervolumeReduction} for objective reduction. DTLZ1 contains $11^k-1$ local Pareto-optimal fronts, and each of them can attract the evolutionary algorithm. Thus, it can be used to test the ability of the algorithm to deal with multiple local attractors. DTLZ2 investigates the performance of the algorithm when dealing with many objective algorithms and the quality of the solutions in regard with proximity to true PF. Finally, DTLZ7 has $2^{m-1}$ disconnected Pareto-optimal regions in the objective space and is used to check the diversity of the solutions and the performance  of the algorithm in disconnected feasible space.

 As suggested~\cite{DebThiLau2005dtlz}, the decision space dimension ($n$) is
$m+4$ for DTLZ1, $m+9$ for DTLZ2 and $m+19$ for DTLZ7. 
In DTLZ problems optimizing a subset of objectives will optimize the rest of the objectives as well. To make the problem more challenging and also to avoid collapsing the PF to one point when projected to $k<m$ objectives, we follow \cite{BroZit2007hypervolumeReduction} and map $x_i$ to $x_i/2 + 0.25$, $i=1,\dotsc,n$, for DTLZ2 and bound $x_i$ within $[0.25, 0.75]$ for DTLZ1, which is also suggested by the authors of BCEMOA~\cite{BatPas2010tec}. Please note that this modification is not needed for DTLZ7 because it does not collapse to a single point. 

\subsection{Machine Decision Maker (MDM)}

We adhere to the MDM framework introduced in \cite{LopKno2015emo} and simulate the DM's preference  with a UF that implicitly expresses which objectives are relevant. 
We define \vec{c} as the ordered index set of relevant objectives such that $i<j \to c_i<c_j$, and consider the following quadratic UFs  that were proposed in experiments on the original BCEMOA~\cite{BatPas2010tec}:
\begin{align}
 \text{UF1}(\vec{f})&= 0.28 f_{c_1}^2 + 0.38f_{c_2}^2 + 0.29 f_{c_1} f_{c_2}  + 0.05 f_{c_1}\label{eq:1}\\
 \text{UF2}(\vec{f}) &= 0.6f_{c_1}^2 + 0.05f_{c_1}f_{c_2}  + 0.23f_{c_1} +  0.38f_{c_2}\label{eq:2}\\
 \text{UF3}(\vec{f})&=  0.44f_{c_1}^2 + 0.14f_{c_2}^2 + 0.09f_{c_1}f_{c_2} + 0.33f_{c_1}  \label{eq:3}
\end{align}
In addition, we consider the following Tchebychef UF:
 \begin{equation}\label{eq:tch}
    U_\text{tch}(\vec{f}) = \max_{i \in \vec{c}}w_i|f_i-f_i^*|
\end{equation}
with $\vec{0}$ as ideal point $f^*$.
The weights $w_i$ of irrelevant objectives ($i \notin \vec{c}$) are set to zero while the weights of relevant objectives are selected by trial and error in a way that the most preferred solution is away from the corner points as far as possible.
As shown in the utilities above, in all cases, the DM only considers
$\Fdm = \{f_i| i \in \vec{c}\}$, while other objectives  are irrelevant. We explain how we selected relevant objectives for each problem in the next section.

Since the DTLZ objectives are to be minimized and to preserve consistency, we assume minimization for utility values as well; i.e. solutions with lower utility value are preferred by the DM. 

\subsection{Selecting relevant objectives}

\begin{figure}[!t]
    \includegraphics[width=0.48\textwidth]{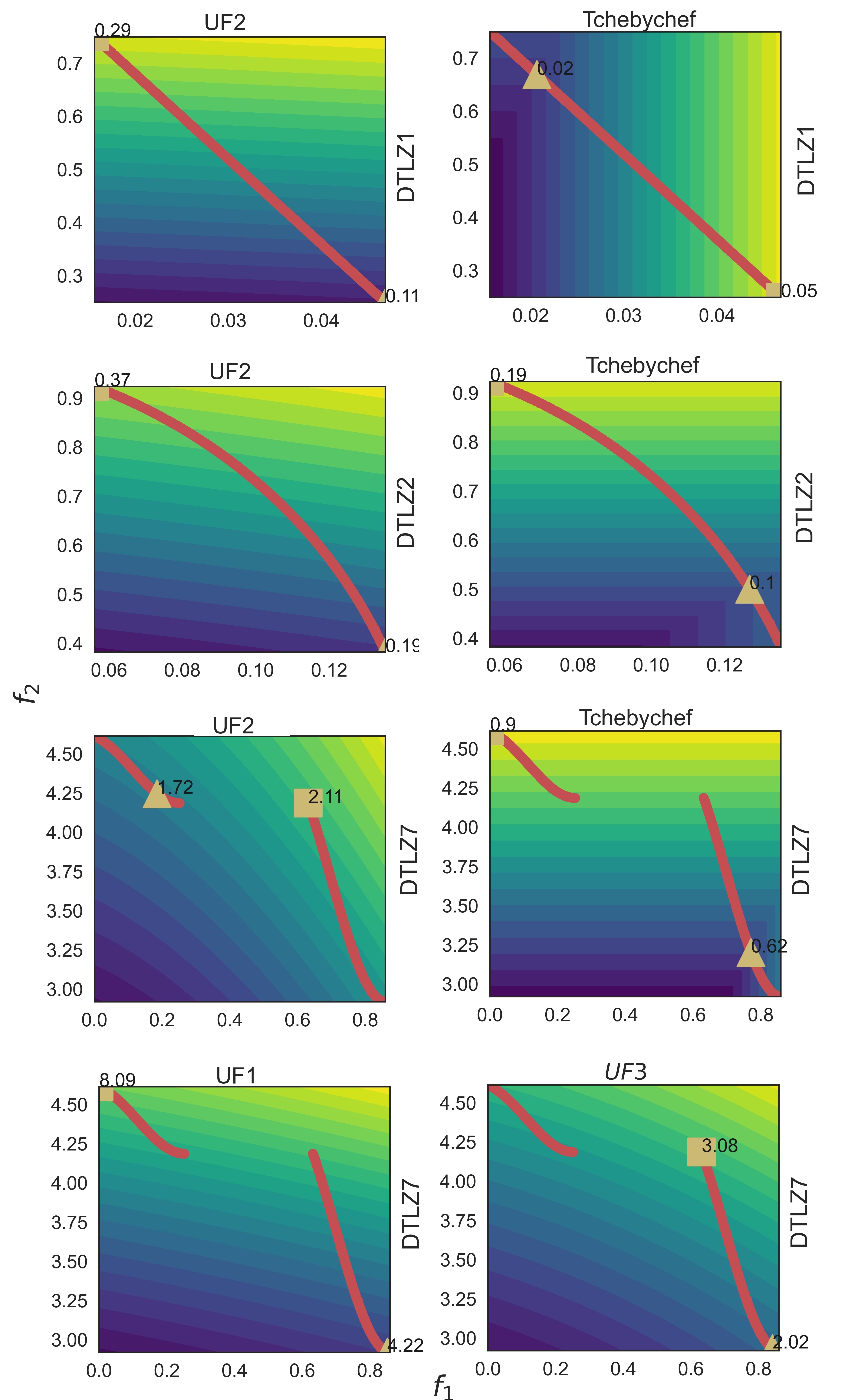}
  \caption{PF of the selected problems with first and fourth objectives being active depicted over the contour lines of the different UFs. The position of the worst and best solutions is marked respectively with a  rectangle and a triangle.}
 \label{fig:PF_contour}
\end{figure}

Projection of the PF on lower dimensions might make it collapse to a single point, for some problems. This is true for DTLZ problems even when the problem is bounded~\cite{BroZit2007hypervolumeReduction}. Thus, the extended problem and its PF should be given a careful examination before proceeding to experiments. For instance, the PF of DTLZ7 would collapse to a point if the first two objectives are active and the rest are inactive. After a careful examination, here the first and fourth objectives are selected as relevant for DTLZ problems ($\vec{c}=\{1,4\}$). In the case of $\rho$MNK problems the first two objectives are selected ($\vec{c}=\{1,2\}$). For  $\rho$MNK problem with four objectives, having $\Fdm = \{f_1, f_2\}$ and given an
 initial $\vec{d} = (0,1,0,1)$, we can see that $f_1$ is a hidden objective
(relevant but not optimized), $f_2$ is both relevant and optimized, $f_3$ is irrelevant and not optimized, and $f_4$ is irrelevant and optimized.
The projection of the PF and the UF contours are depicted in Figure~\ref{fig:PF_contour} for DTLZ problems and selected relevant objectives. The location of the most preferred solution is indicated by a triangle and the worst solution with a rectangle. The graphs for UF1 and UF3 are identical to UF2 on DTLZ1 and hence omitted to save space. The same is true for DTLZ2.

\subsection{Evaluation of the results}
The experiments are performed in three different modes to enable the assessment of the algorithms:
\begin{enumerate}
\item \textit{Golden} mode: No interaction is done in this mode and the algorithm directly accesses the true UF of the DM instead of learning a UF. Moreover, only relevant objectives
  are optimised from the start to end. This is the ideal scenario.

\item \textit{Only Learning}  mode: This mode corresponds to the original BCEMOA without any detection of hidden objectives. The algorithm does not have access 
 to the DM's UF and
  instead a UF is learned from pairwise comparisons provided by the MDM at
  each  interaction; i.e. at each interaction the MDM uses true UF to rank solutions. Predictions from the learned utility are used to rank
  non-dominated solutions, replacing the crowding distance in NSGA-II.  The
  algorithm still uses non-dominated sorting as the first criteria to rank
  solutions. Both non-dominated sorting and the learned utility only consider the set of active objectives $\hat{\vec{f}}(x)$.  The set of
  active objectives never changes, that is, $\vec{d}$ remains constant
  throughout the run. 
  
\item \textit{Learning $+$ detection} mode: This is our proposed BCEMOA-HD
  that performs detection of hidden objectives and is able to modify the set of
  active objectives. Within this mode, we test $4$ variants of the HD method: \khd, \thd, \khdr, \thdr. Similar to the \textit{Only Learning} mode, the optimization algorithm relies
  on non-dominated sorting and an UF that is learned based on $\hat{\vec{f}}(x)$, and not the DM's true utility. However, in this mode, $\vec{d}$ and subsequently $\vec{\hat{f}}$ may change
  after each DM interaction with the ultimate goal of converging to the objectives that are actually relevant for the DM ($\Fdm$).
\end{enumerate}

Having these three modes makes it possible to evaluate the performance of the proposed method compared to the original BCEMOA and to the best solution that is achieved under an ideal scenario in \textit{Golden} mode. The criteria for evaluation of the performance is the true utility value of the final solution returned by the algorithm in different modes. We also make record of the active objectives after each interaction for \textit{Learning $+$ detection} mode to investigate how good the proposed method performs in detecting and diverging towards the true set of objectives.

\subsubsection{Parameter settings}
In all evaluation modes, BCEMOA uses the parameter settings proposed in the original
paper~\cite{BatPas2010tec}, including the parameters of the SVM learning model. In particular, it uses a population size
of $100$ and creates $100$ new solutions at each generation, The total number of generations is 500 and $\Nexa=5$ solutions are shown to the DM
at each interaction.

Within BCEMOA, NSGA-II  runs for
$gen_1=200$ generations before the first interaction and there are $gen_i=30$ generations between subsequent interactions. The total number of generations after the last interaction is calculated by $500-gen_1- gen_i(N_\text{int}-1)$. Thus, changing the number of interactions ($N_\text{int}$) would not alter total number of generations. 
We run experiments with $1$, $3$  and $6$ interactions for DTLZ problems. For $\rho$MNK problems, we only consider $6$ interactions and, instead, we investigate the effect of different levels of correlation ($\rho$) and ruggedness($K$).
Each test was repeated 40 times with different random seed.

\subsubsection{Implementations}
The algorithms, machine DM and $\rho$MNK problems are implemented in Python
3.7.6. The implementations of NSGA-II  within BCEMOA and DTLZ benchmarks are provided by the Pygmo library 2.16.0 \cite{BisIzzYam2010:pagmo-arxiv}, the uni-variate feature-selection and RFE implementations are based on Scikit-learn
0.23.1 (\url{http://scikit-learn.org/}).  

\section{Experimental Results \& Discussion} \label{results}
We have considered $3$ DTLZ benchmark problems and each of them are considered with $3$ different values for number of interactions. The $\rho$MNK problems are considered with $6$ values of $\rho$ and $4$ values of $K$. Considering $3$ different number of objectives, in total we have investigated $99$ problems, each of which are considered with $4$ utility functions. Each experiment is repeated $40$ times with different random seeds and the average values are reported.
 Comparison of different modes is done with regard to the utility value of the final solution returned by the algorithm averaged over $40$ runs.
 In this section we focus on the most important findings.
Figures that do not include key findings are not presented to save space. However, the complete set of results and figures can be found in the Appendix. 

\subsection{DTLZ problems with fixed number of objectives}
\begin{figure}[!t]
    \includegraphics[width=0.48\textwidth]{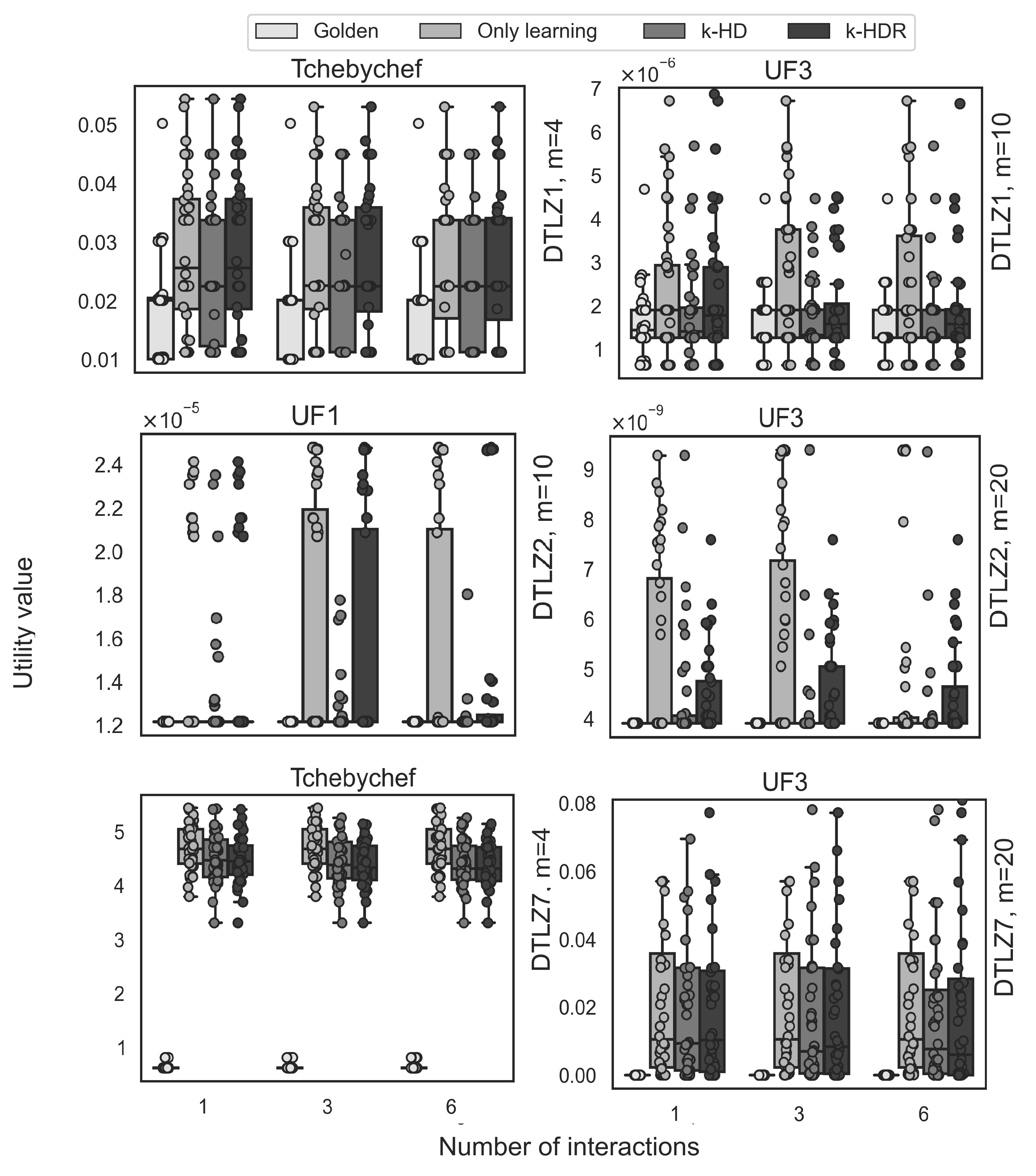}
 \caption{Comparison of the performance of different modes for DTLZ problems. Vertical axis indicated the utility value. Number of active objectives are fixed. Horizontal axis indicates number of interactions. Lower utility values are preferred.}
 \label{fig:127fixed}
\end{figure}

The results of experiments on DTLZ problems with fixed number of active objectives is illustrated in Fig.~\ref{fig:127fixed}.
When using online detection of hidden objectives with DTLZ1 problem with $m=4$ or $m=20$ and fixed number of active objectives, almost no improvement  is observed in terms of the utility value compared to \textit{Only learning}. However, when $m=10$ the performance of \khd and \khdr is significantly better than \textit{Only learning} and almost as good as \textit{Golden} mode. The only exception is when UF2 is used. 

For DTLZ2, improvements in the performance can be seen when detection of hidden objectives is active in the case of $m=10$ and $m=20$ when utility functions UF1 and UF2 are used. Another important observation is better performance of \khdr with more interactions, although it fails to get as good as \khd.

For DTLZ7, there are slight improvements when detection methods (\khd, \khdr) are used. Complete list of figures can be found in the Appendix. In general, it is observed that the proposed method can significantly improve the utility value of the final solutions. Although in some cases it fails to do so, the utility value is not deteriorated by the method.

\subsection{DTLZ problems with variable number of objectives}
\begin{figure}[!t]
    \includegraphics[width=0.48\textwidth]{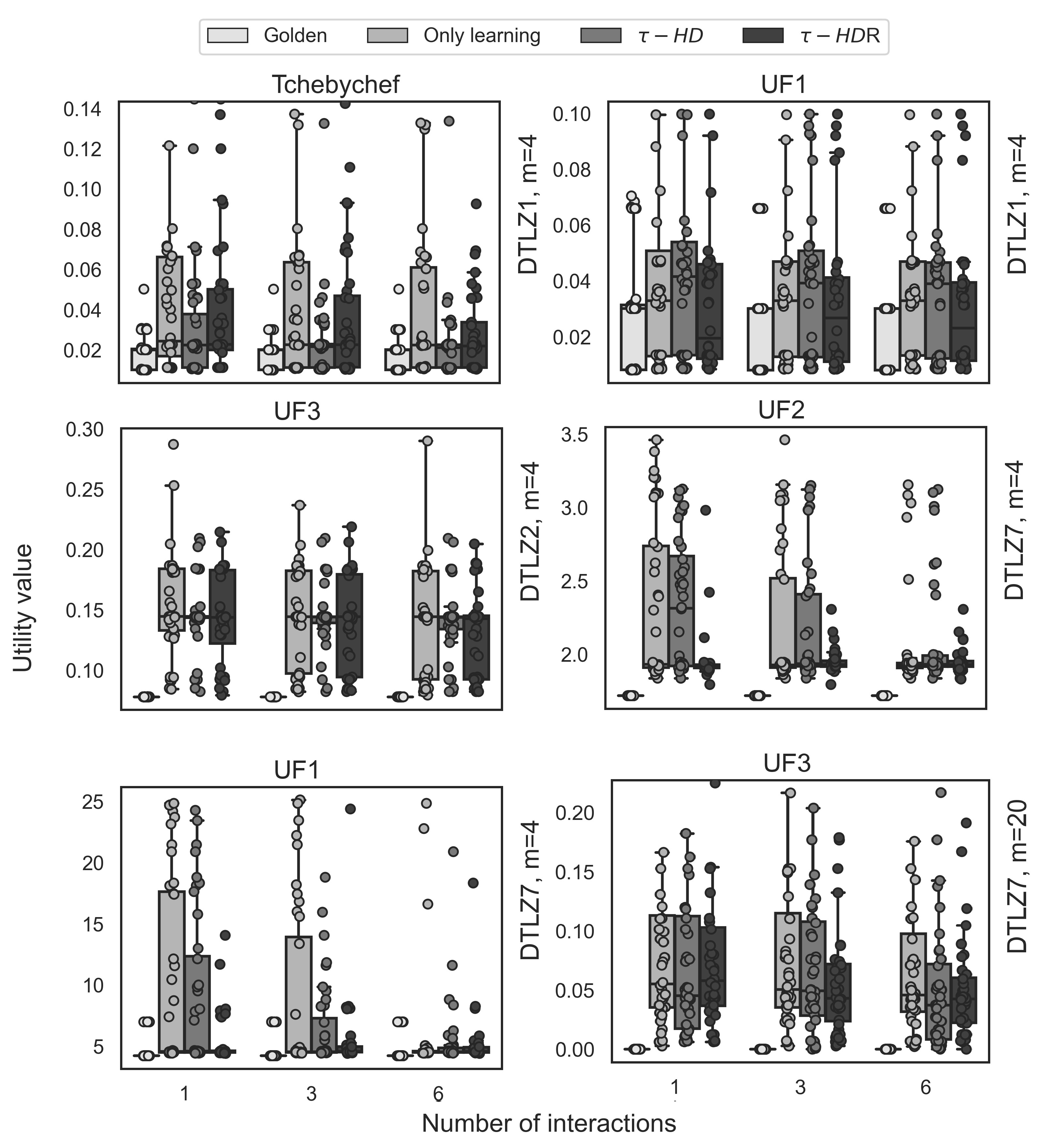}
 \caption{Comparison of the performance of different modes for DTLZ problems. The number of active objectives is not fixed and detection mode is used as a mean of objective reduction technique. Vertical axis indicated the utility value. Horizontal axis indicates number of interactions. Lower utility values are preferred.
 }
 \label{fig:127reduction}
\end{figure}

In this set of experiments, the effectiveness of the proposed detection method is investigated with regard to objective reduction capabilities and thus the number of active objectives is not fixed. Furthermore, the experiments start with all the objectives being active. The key results of these experiments are illustrated in Figure~\ref{fig:127reduction}. 
For DTLZ1 with $m=4$, the \thd and \thdr perform better than the \textit{only learning} mode on Tchebychef UF, while for other UFs they have almost the same performance. With $m=10$ and $m=20$, \thd and \thdr performs as good the \textit{Golden} mode while $\tau$ is slightly outperformed by \thd. Results for DTLZ2, are identical to those of DTLZ1; i.e. \thd and \thdr still doing as good as \textit{Golden} mode and outperform \textit{Only learning} mode with $m=10$ and $m=20$. For $m=4$ the performance of \thd and \thdr outperforms \textit{only learning} when UF3 is used, but cannot perform as good as \textit{Golden} mode. 

For DTLZ7 with $m=4$ all algorithms perform as good as \textit{Golden} mode which suggests the problem may be too easy for the algorithms and optimizing any subset of objectives may lead to optimization of other objectives as well~\cite{BroZit2007hypervolumeReduction}. However, with higher number of objectives, they fail to perform as good as the \textit{Golden} mode and \thd and \thdr outperform \textit{Only learning}.

\subsection{$\rho$MNK  problems}
\begin{figure}[]
    \includegraphics[width=0.48\textwidth]{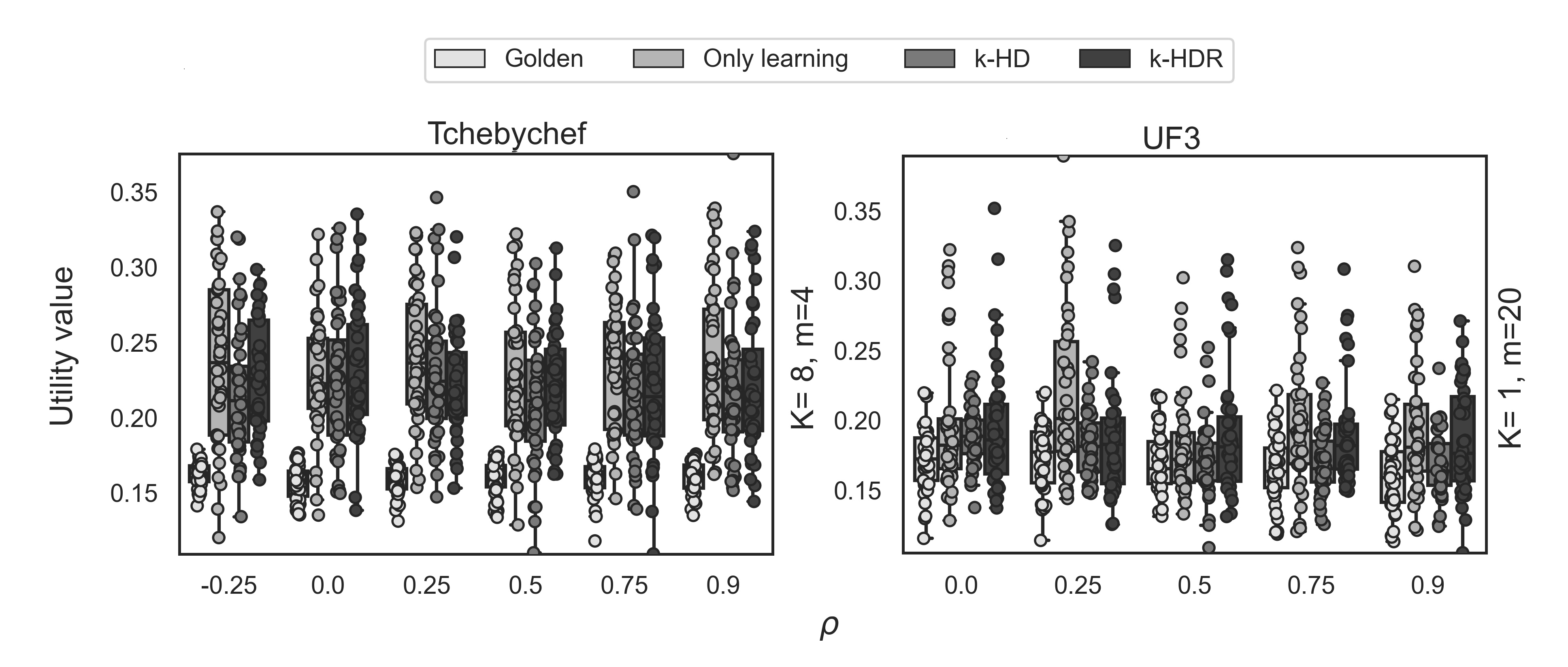}
 \caption{Comparison of the performance of different modes for $\rho$MNK problems. Number of active objectives is fixed. Vertical axis indicated the utility value. Horizontal axis indicates the $\rho$ values. Lower utility values are preferred.}
 \label{fig:rmnkFixed}
\end{figure}

\begin{figure}[!t]
    \includegraphics[width=0.48\textwidth]{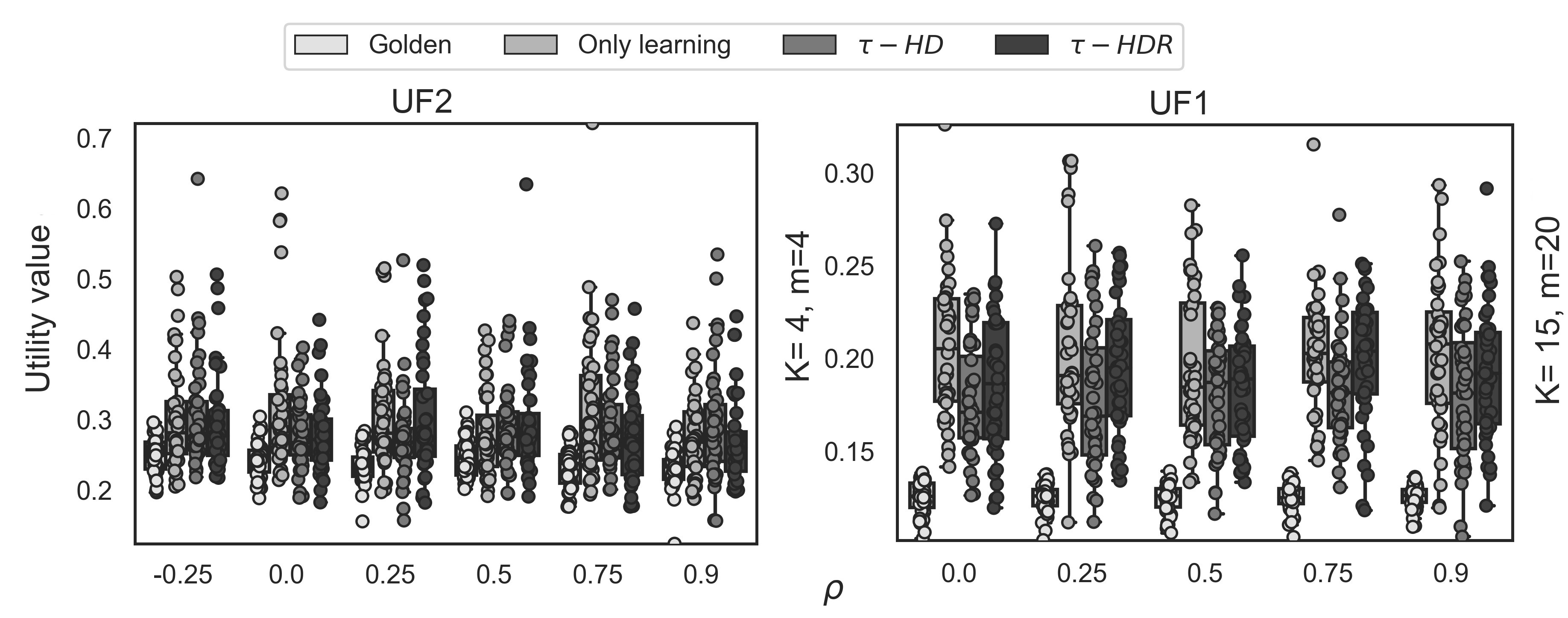}
 \caption{Comparison of the performance of different modes for $\rho$MNK problems. The number of active objectives is not fixed and detection mode is used as an objective reduction technique. Vertical axis indicated the utility value. Horizontal axis indicates the $\rho$ values. Lower utility values are preferred.}
 \label{fig:rmnkreduction}
\end{figure}

The results of experiments on $\rho$MNK problems with fixed number of objectives and variable number of objectives are depicted in Figures~\ref{fig:rmnkFixed} and \ref{fig:rmnkreduction}, respectively. In all of these experiments, detection method fails to make significant improvements in terms of the utility value compared to \textit{Only learning} mode. However, there are significant savings in terms of computational cost when \thd and \thdr are used. In the next Section, this aspect will be discussed in full.

\subsection{Sensitivity analysis and Discussions}\label{sensitivity}
\subsubsection{Power of Detection of True Set of Objectives}

\begin{figure}[!t]
    \includegraphics[width=0.49\textwidth]{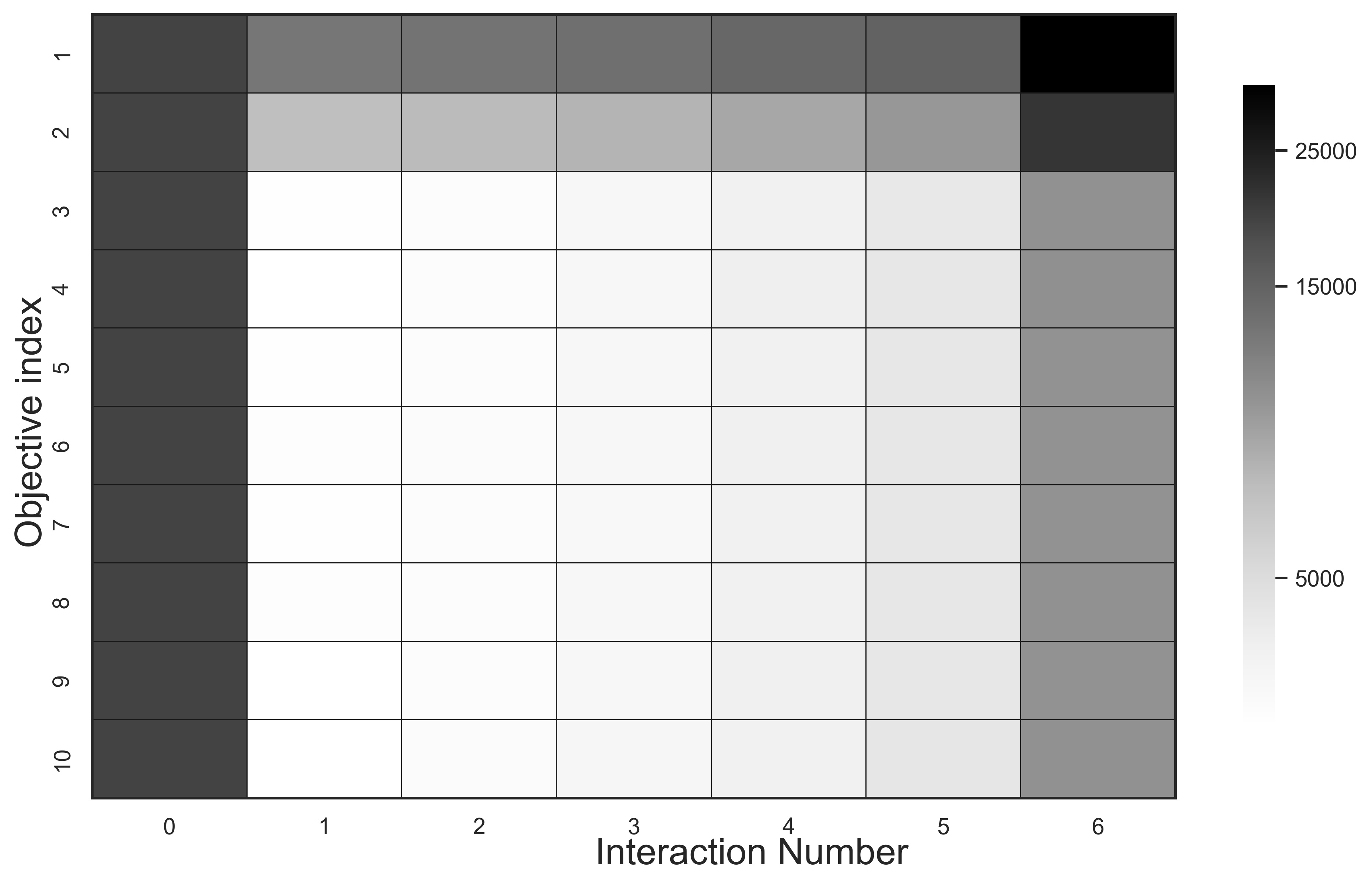}
 \caption{Divergence of the algorithm towards true objective set through interactions on problem $\rho$MNK with $10$ objectives ($c=\{1,2\}$).  The x-axis indicates the interaction number and the y-Axis indicates the index of the objective functions. The darker color, the higher is the number of the times the objective has been activated. The true set of objectives is found after the first interaction in the majority of the experiments.}
 \label{fig:a_heatmap}
\end{figure}

In terms of the power of the detection, a heatmap plot is provided in Figure~\ref{fig:a_heatmap}. The plot illustrates the frequency of the times the true set of objectives are activated by the detection of hidden objectives across all experiments on $\rho$MNK problems with $10$ objectives when \thd is on. The x-axis indicates the number of interaction within a single run and the y-axis (rows) pertain to different objectives numbered from $1$ to $10$. Interaction $0$ refers to the state of the algorithm before the first interaction, when all objectives are active. At the first interaction all objectives are active. Then immediately after the first interaction most of the objectives are dropped and first and second objectives that are relevant and included in the utility functions, are kept active. It can be easily verified that the \thd converges fast towards the true set of objectives. Another surprising observation is that after the $6^{th}$ interaction, almost all objectives become active, although to a lesser degree comparing with the relevant ones. This observation can be justified as follows: When detection of hidden objectives is enabled, the selected objectives might be optimized to near-optimal values and then become fixed in the next generations. The uni-variate correlation of these objectives having fixed values with ranking vector would diminish to almost $0$. As a result, the objectives would be identified as irrelevant and substituted with new ones whether or not they are important to the DM.

\subsubsection{Analysis of Performance and Threshold $\tau$}
\begin{figure*}[!t]
    \includegraphics[width=0.98\textwidth]{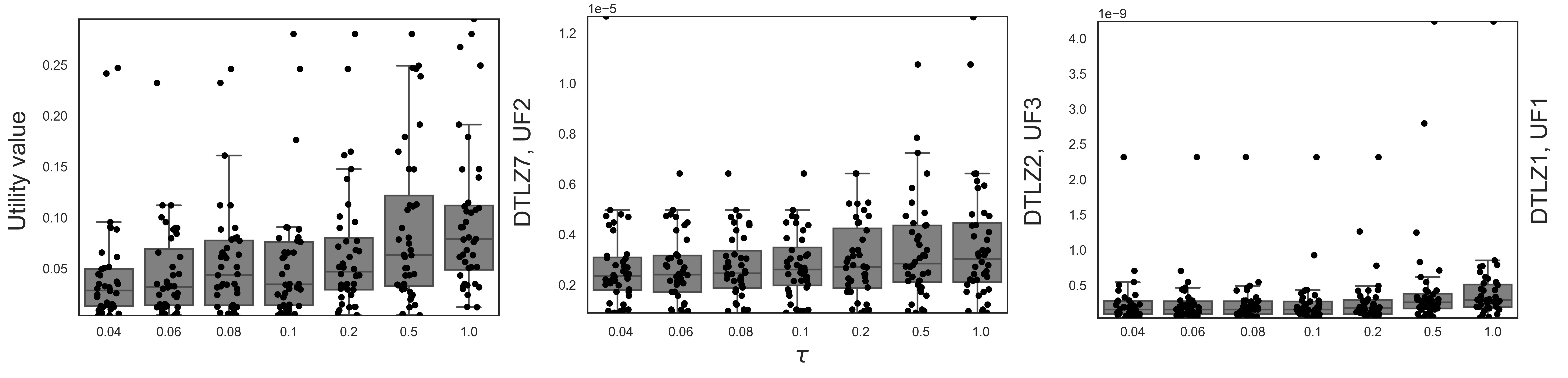}
 \caption{Performance analysis of \thdr with different values of $\tau$ for DTLZ problems with $m=20$. Number of interaction in all runs is fixed to 6.}
 \label{fig:DTLZ_tau_analysis}
\end{figure*}

As an important parameter of \thd  ~and \thdr, $\tau$ has a direct effect on the number of active objectives and those that are removed. Thus, careful examination should be given in determining this parameter. 
To inspect the effect of parameter $\tau$, which indirectly controls the number of active objectives, the DTLZ problems with $m=20$ objectives are solved for different values of $\tau$. The results for other problems are similar and hence not discussed here. Setting $m=20$ would provide for a better illustration of the efficiency of the proposed method in reducing the computational requirements and objective evaluations. When $\tau=1$, all objectives have $p$-value less than the threshold and, thus, all of them are active; this means no objective reduction is performed and the mode is identical to \textit{Only learning}. The results in Figure~\ref{fig:DTLZ_tau_analysis} show that the performance of the \thd improves with lower values of $\tau$ on DTLZ problems. 
On the other hand, reducing the $\tau$ value would reduce the number of active objectives. As explained earlier, in the case where there are fewer than $2$ objectives that pass the selection criteria, the two objectives with least $p$-values are selected.

\begin{figure}[!t]
    \includegraphics[width=0.49\textwidth]{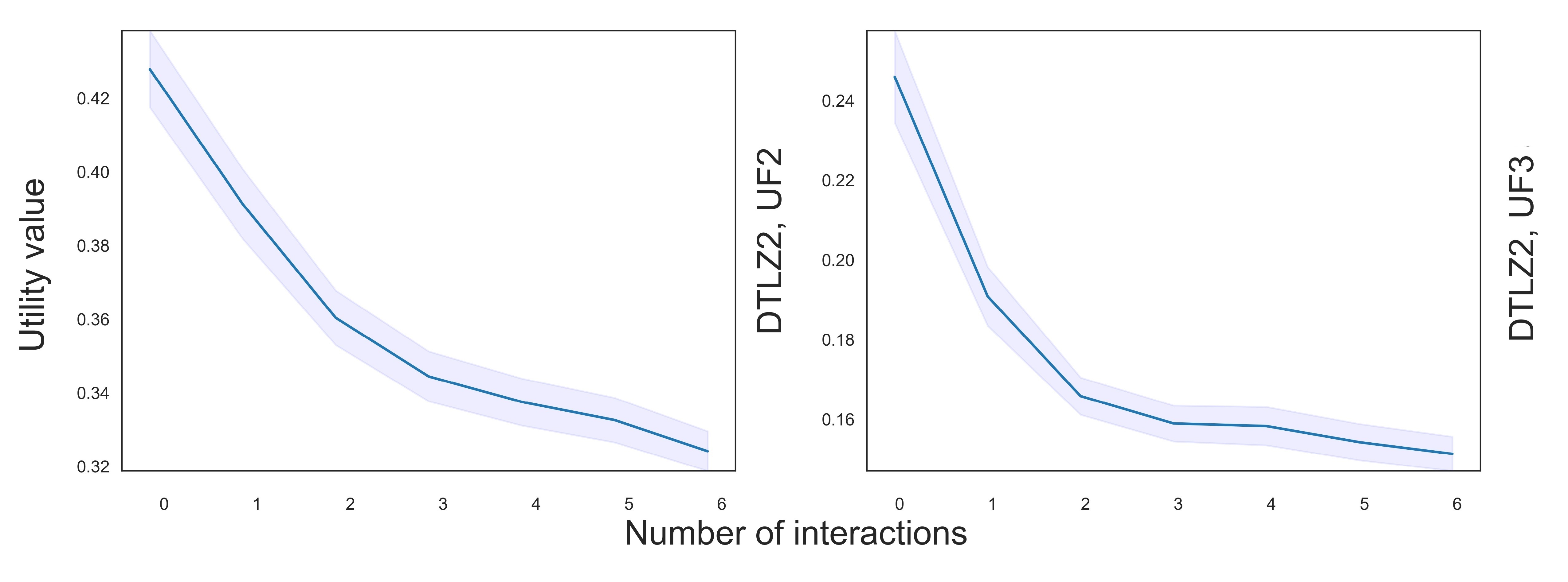}
 \caption{Utility of the best-so-far solution within a single run after each interaction on DTLZ2 problem when  \thd detection mode is on. Results are averaged over 40 runs. The results for other test problems are identical.}
 \label{fig:trend_dtlz}
\end{figure}

Figure \ref{fig:trend_dtlz} illustrates the change in the utility value of the best solution gained after each interaction in a single run of the algorithm averaged over 40 runs. All changes in the objectives results in improvement of the utility value after each interaction. The changes become smaller as the algorithm converges towards a good solution.

\begin{figure}[!t]
    \includegraphics[width=0.49\textwidth]{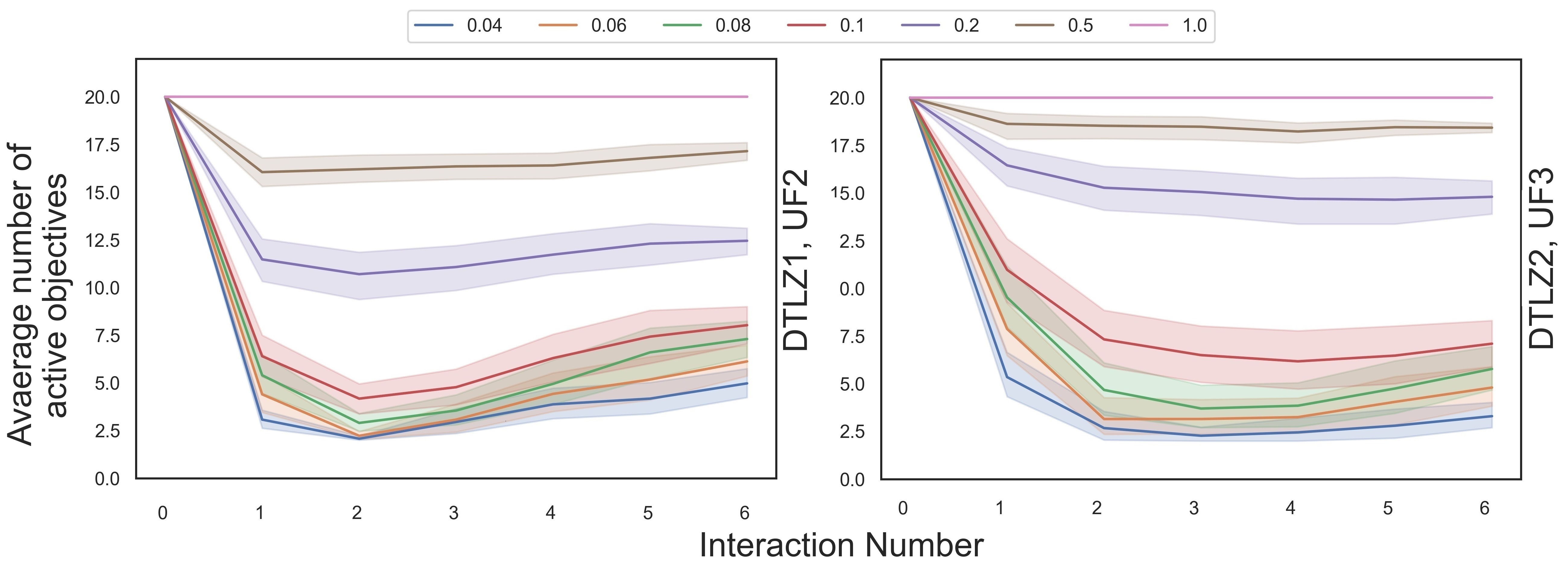}
 \caption{Number of active objectives within a single run after each interaction for different values of $\tau$ in \thdr mode on DTLZ1 and DTLZ2 test problems. The x-axis indicates the number of the interaction and the y-axis illustrates the mean value of the total number of active objectives averaged over 40 runs. The shaded areas show the 95\% confidence interval around the mean. The results for other test problems are identical.}
 \label{fig:a_analysis_DTLZ}
\end{figure}

The changes in the number of active objectives after each interaction, averaged over $40$ runs, are depicted in Figure \ref{fig:a_analysis_DTLZ} where the shaded areas show the 95\% confidence interval around the mean. It can clearly be verified that after the first interaction, the number of active objectives experiences a steep decrease. As expected, when $\tau=1$, no objective reduction is performed.


Lower number of active objectives translates to lower computational costs. However, another important criteria to consider is the ratio of evaluations of relevant objectives to those of irrelevant ones. 

\thd effectively reduces the number of objective evaluations and most of objective evaluations pertain to relevant objectives. As defined, whenever we refer to objective evaluations, we measure the total number of individual objective evaluations. Thus, the number of objective evaluations depend on the number of active objectives.
For instance, when $\tau=1$ (equivalent to \textit{Only learning} mode), an experiment on DTLZ1 with UF3 and $\tau=0.2$,
uses $600,000$ objective evaluations after the first interaction (before the first interaction everything is identical, thus no comparison is made) and only $10\%$ of these evaluations pertain to relevant objectives. However, when \thd is used, only $45,000$ objective evaluations are done of which $30,000$ $(67\%)$ are dedicated to relevant objectives. In general, objective evaluations are reduced by up to $80\%$ compared to \textit{Only learning} when \thd or \thdr is used.

These results prove the effectiveness of the \thdr in terms of both performance with regard to the quality of the final results and computational efficiency.

\section{Conclusion and future work} \label{conclusion}
This study has considered interactive multi-objective problems where only an unknown subset of all the defined objectives are of relevance to the DM. In this context, we provided formal definitions of irrelevant, hidden and active objectives that complement the definition of redundant objectives already studied in the literature.
Consequently, simulation of irrelevant, hidden, and active objectives are discussed and an efficient method is proposed to deal with these objectives. 

Two variants of the method with variable and fixed number of active objectives were studied. The results show that the variant with variable number of objectives can be considered for dimension reduction purposes, reducing efficiently the number of active objectives even after the first interaction; this eliminates unnecessary evaluations of irrelevant objectives and improves computational time and quality of the final solution returned by the algorithm. The variant with fixed number of active objectives also manages to detect and switch to relevant objectives. We also explored the application of recursive feature selection to investigate if further improvements could be achieved. However, the results indicate that there is no gain in using this method over the uni-variate feature selection. Comparing the results achieved for different test problem suits, it can be observed that improvement in the utility value of the final solutions are more significant for DTLZ problems. However, achievements with regard to saving in objective evaluations and computational costs

is similar for both test problems. The benefits of this achievement is highlighted where expensive objective evaluations are avoided. The effect of correlation among objectives was investigated and proved to be insignificant in tests with $\rho$MNK problems.

A sensitivity analysis was performed to scrutinize different aspects of the problem, effect of key parameters of the algorithms on final results, and to provide further insight into the proposed method. We showed experimentally that the value of $\tau$ affects the number of active objectives and can be used as a tool to control this aspect.
The proposed feature selection method only considers linear relations between objectives and DM's rankings. Considering nonlinear regression in feature selection would be a subject worth to study, although it adds to the computational effort. In the case of BCEMOA, one can use the SVM's support vector in order to identify most important objectives. However, the intention of this study was to target ranking-based algorithms in general rather than a specific one. 
We considered 4 different UFs to capture MDMs with different preferences. Future study could expand on this and consider a range of further UFs, such as a Sigmoid UF and others \cite{ShaLopKno2021gecco}. 
For some experiments, we observed that once the relevant objectives
are optimized, the proposed methods try to explore optimizing
other inactive objectives in the hope to find an even better utility. In future studies, it would be desirable to introduce a stopping mechanism  to avoid such excess search and thus increase computational efficiency/reduce resource usage.

{\footnotesize
  \bibliographystyle{IEEEtranN}
  \bibliography{bib/abbrev,bib/journals,bib/authors,bib/biblio,bib/crossref,bib/my_bib}
}
\clearpage

\end{document}


\onecolumn
\bstctlcite{BSTcontrol}
\title{Detection of Hidden objectives and Interactive Objective Reduction\\Supplementary Material}
\author[1]{Seyed Mahdi Shavarani}
\author[1]{Manuel López-Ibáñez}
\author[1]{Richard Allmendinger}
\affil[1]{Alliance~Manchester~Business~School,~University~of~Manchester, Manchester,~M13~9PL,~UK \authorcr Email: {\tt \{seyedmahdi.shavarani,manuel.lopez-ibanez, richard.allmendinger\}@manchester.ac.uk}\vspace{1.5ex}}


\maketitle

This supplementary material provides all the Figures that were discarded from the main text for brevity. 

\thispagestyle{empty}

\listoffigures

\clearpage
\pagenumbering{arabic}

\begin{figure*}[!t]
    \includegraphics[width=0.8\textwidth]{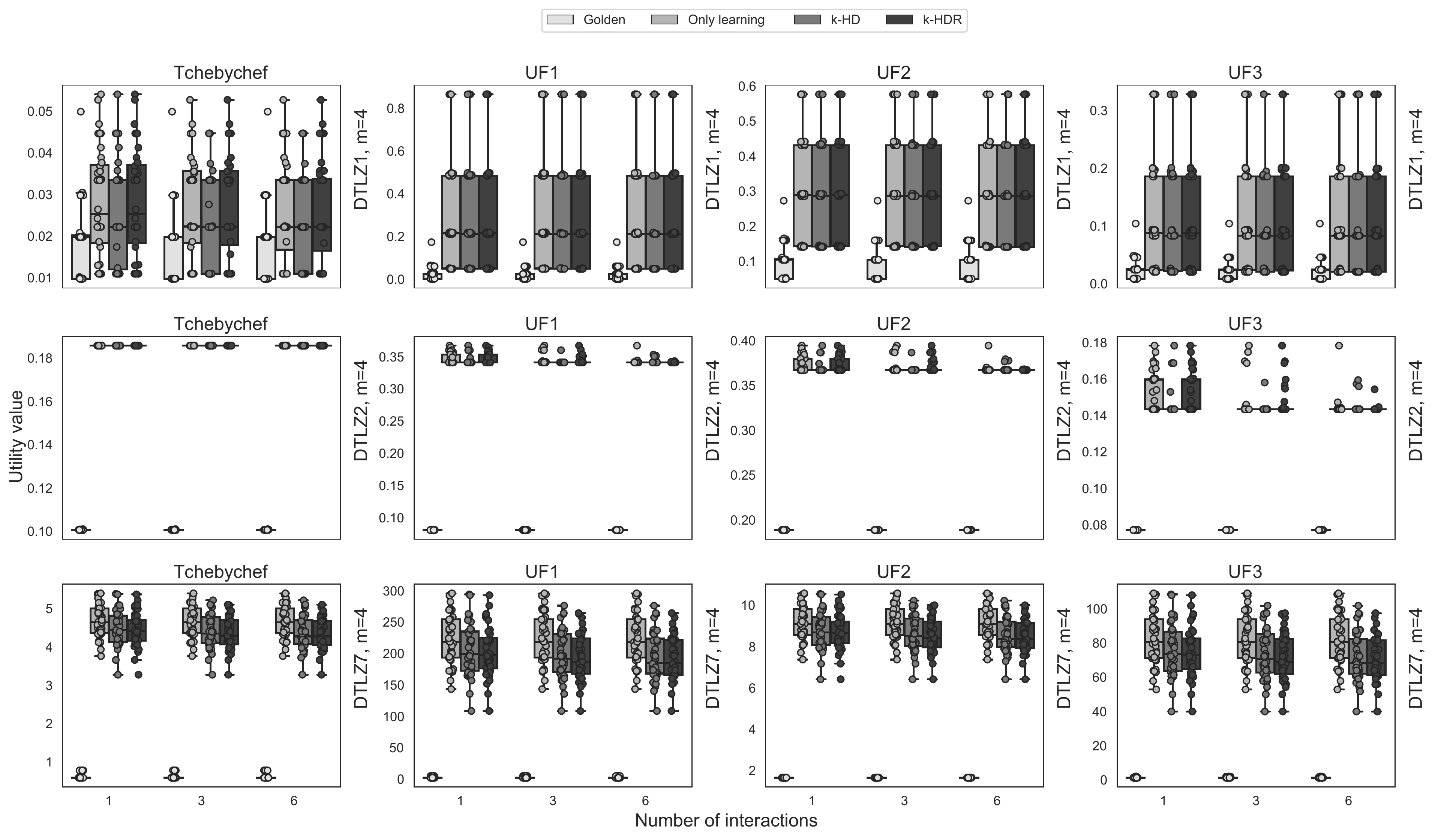}
 \caption{Comparison of the performance of different modes for DTLZ problems with $m=4$. The number of active objectives are fixed.
 }
 \label{fig:127fixed4}
\end{figure*}

\begin{figure*}[!t]
    \includegraphics[width=0.8\textwidth]{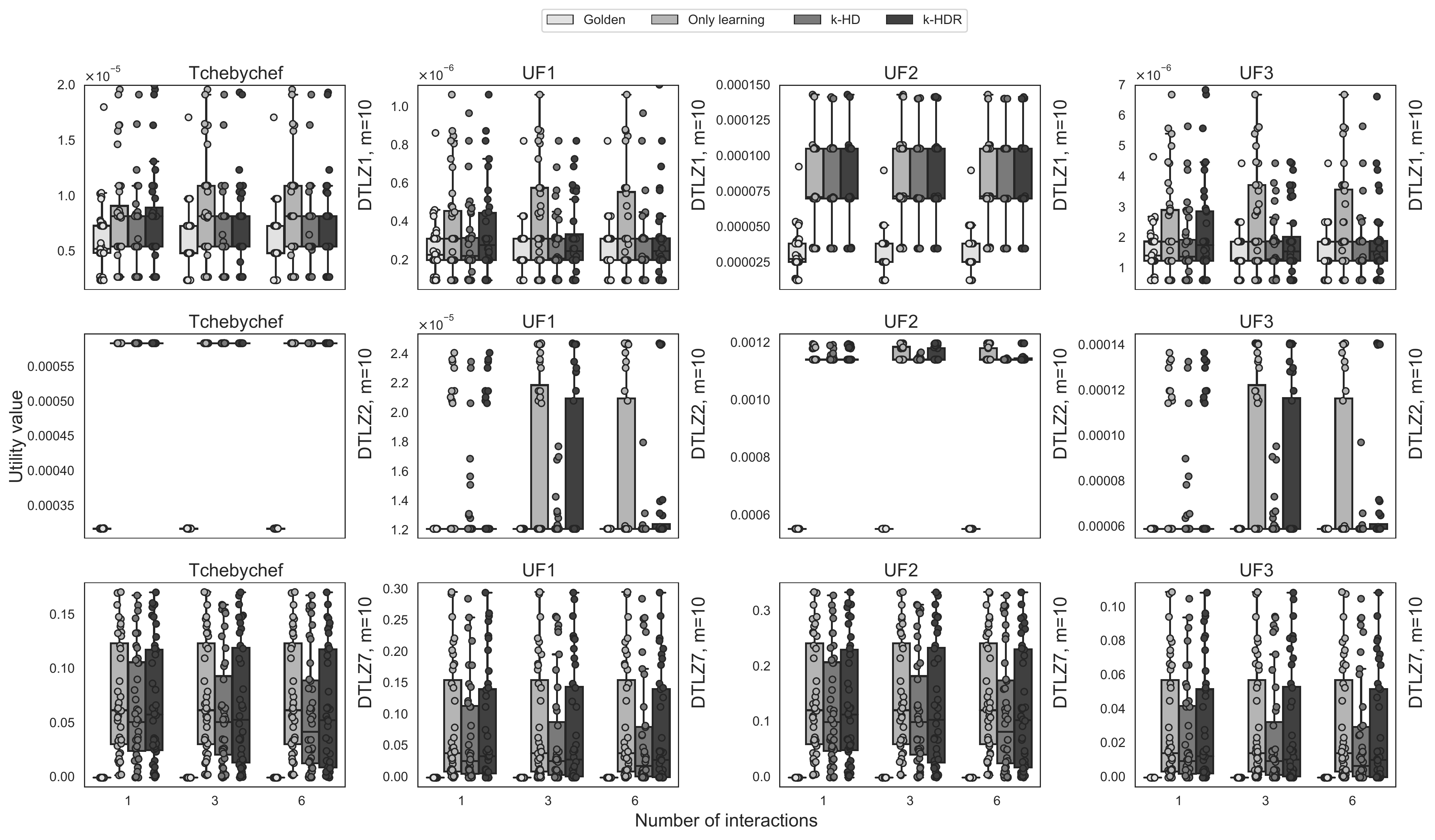}
 \caption{Comparison of the performance of different modes for DTLZ problems with $m=10$. The number of active objectives are fixed.
 }
 \label{fig:127fixed10}
\end{figure*}

\begin{figure*}[!t]
    \includegraphics[width=0.8\textwidth]{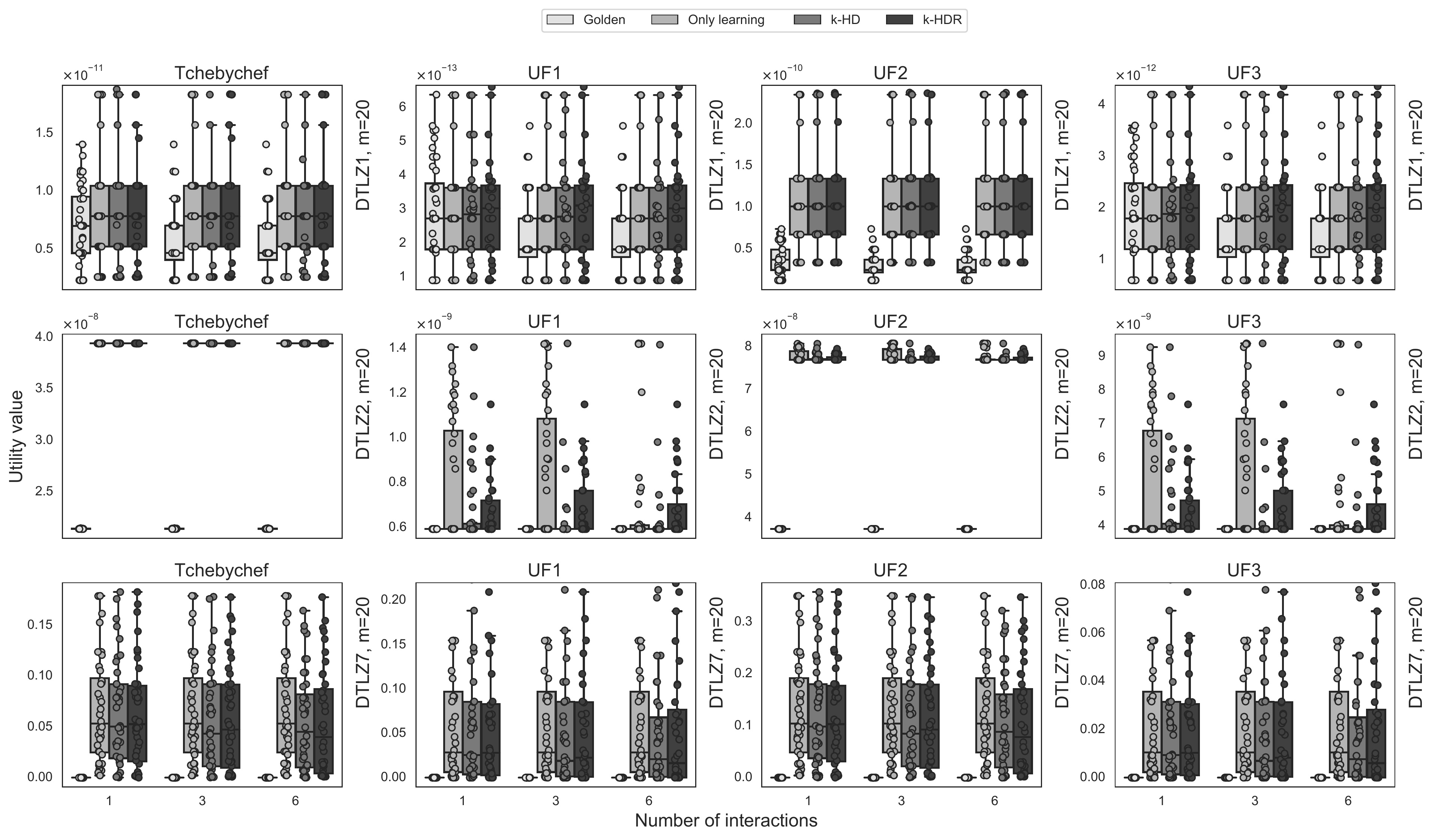}
 \caption{Comparison of the performance of different modes for DTLZ problems with $m=20$. The number of active objectives are fixed.
 }
 \label{fig:127fixed20}
\end{figure*}

\begin{figure*}[!t]
    \includegraphics[width=0.8\textwidth]{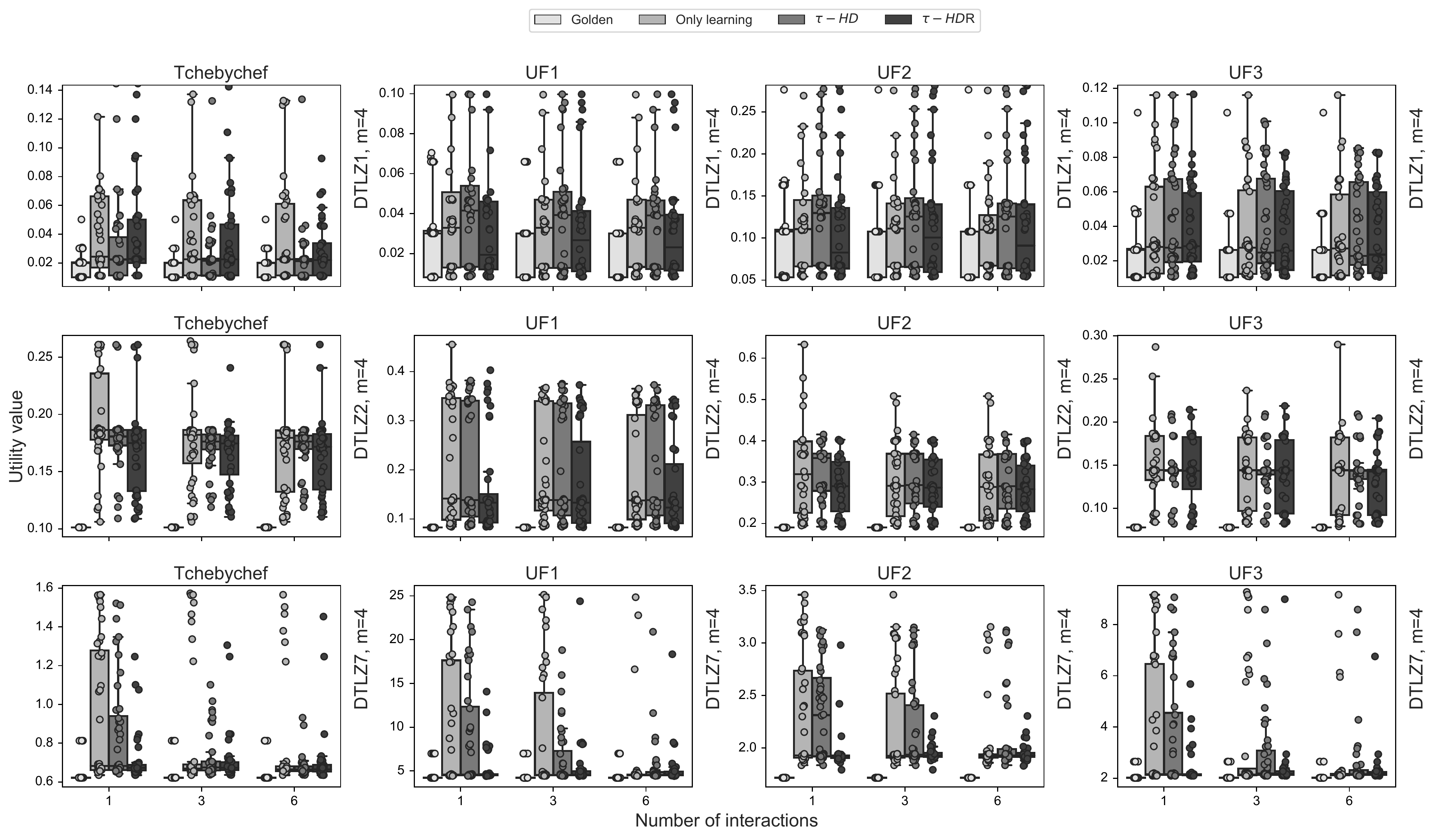}
 \caption{Comparison of the performance of different modes for DTLZ problems with $m=4$. The number of active objectives is not fixed and detection mode is used as a mean of objective reduction technique.
 }
 \label{fig:127reduction4}
\end{figure*}

\begin{figure*}[!t]
    \includegraphics[width=0.8\textwidth]{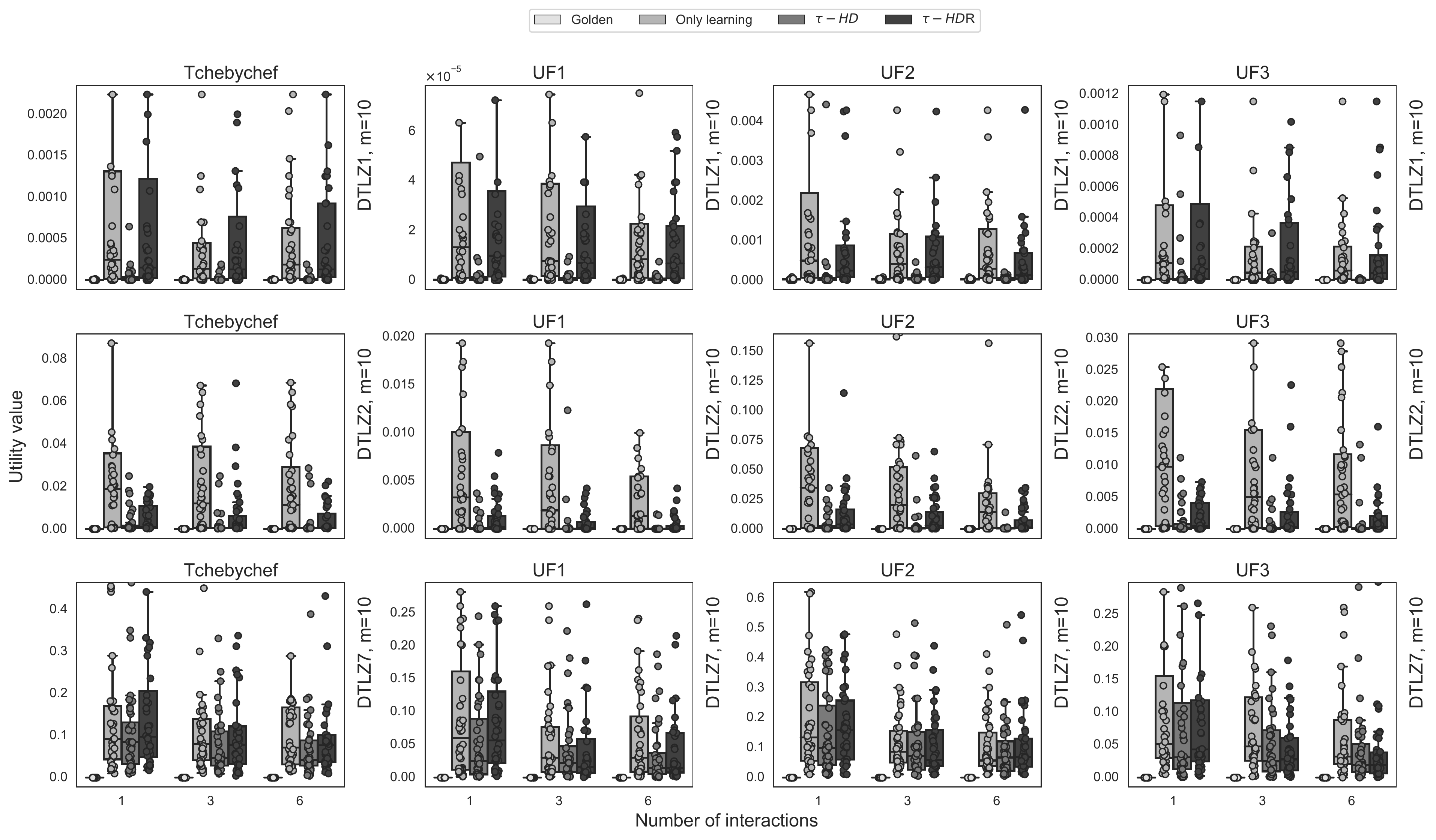}
 \caption{Comparison of the performance of different modes for DTLZ problems with $m=10$. The number of active objectives is not fixed and detection mode is used as a mean of objective reduction technique.
 }
 \label{fig:127reduction10}
\end{figure*}

\begin{figure*}[!t]
    \includegraphics[width=0.8\textwidth]{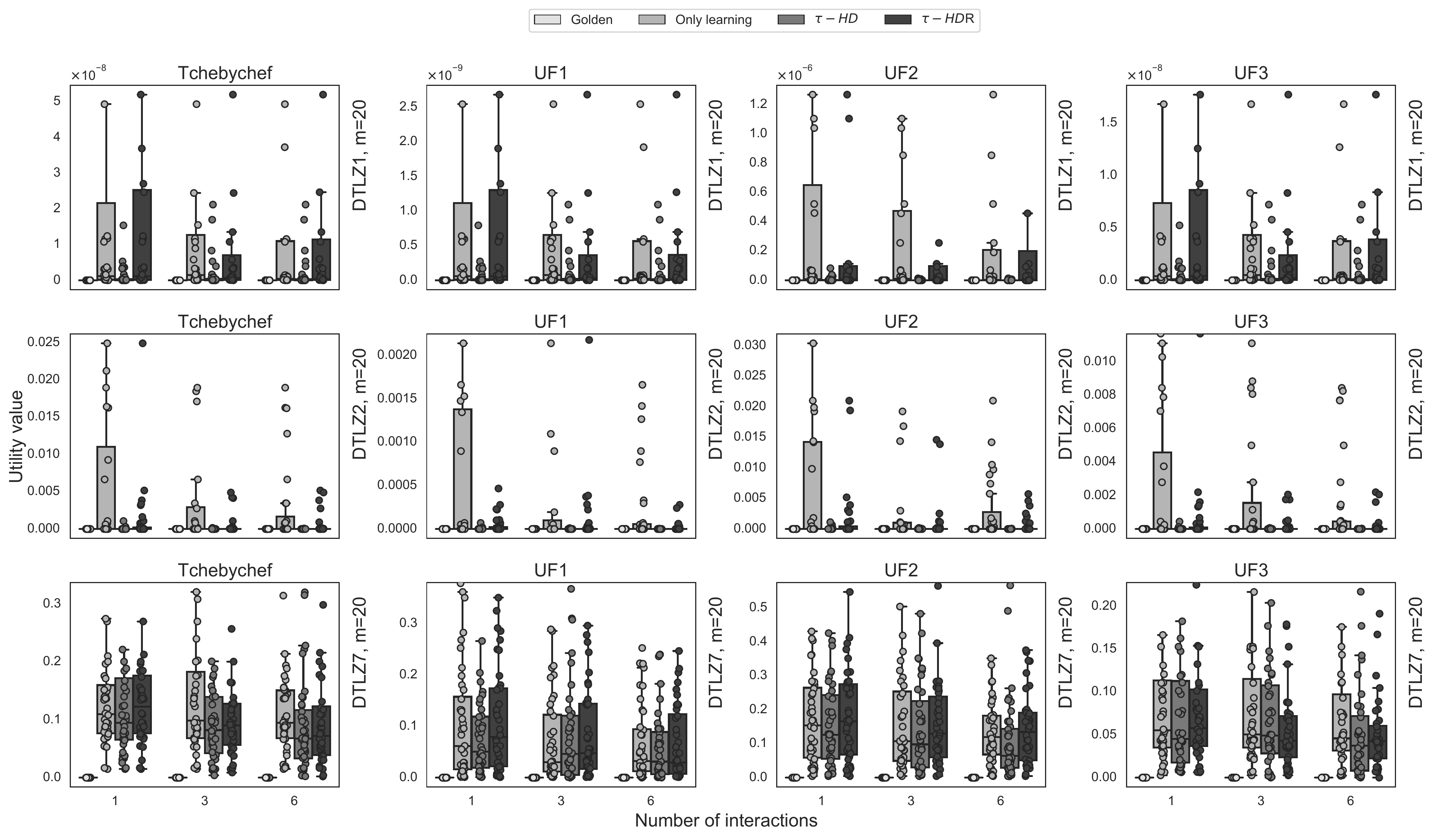}
 \caption{Comparison of the performance of different modes for DTLZ problems with $m=20$. The number of active objectives is not fixed and detection mode is used as a mean of objective reduction technique.
 }
 \label{fig:127reduction20}
\end{figure*}

\begin{figure*}[!t]
    \includegraphics[width=0.8\textwidth]{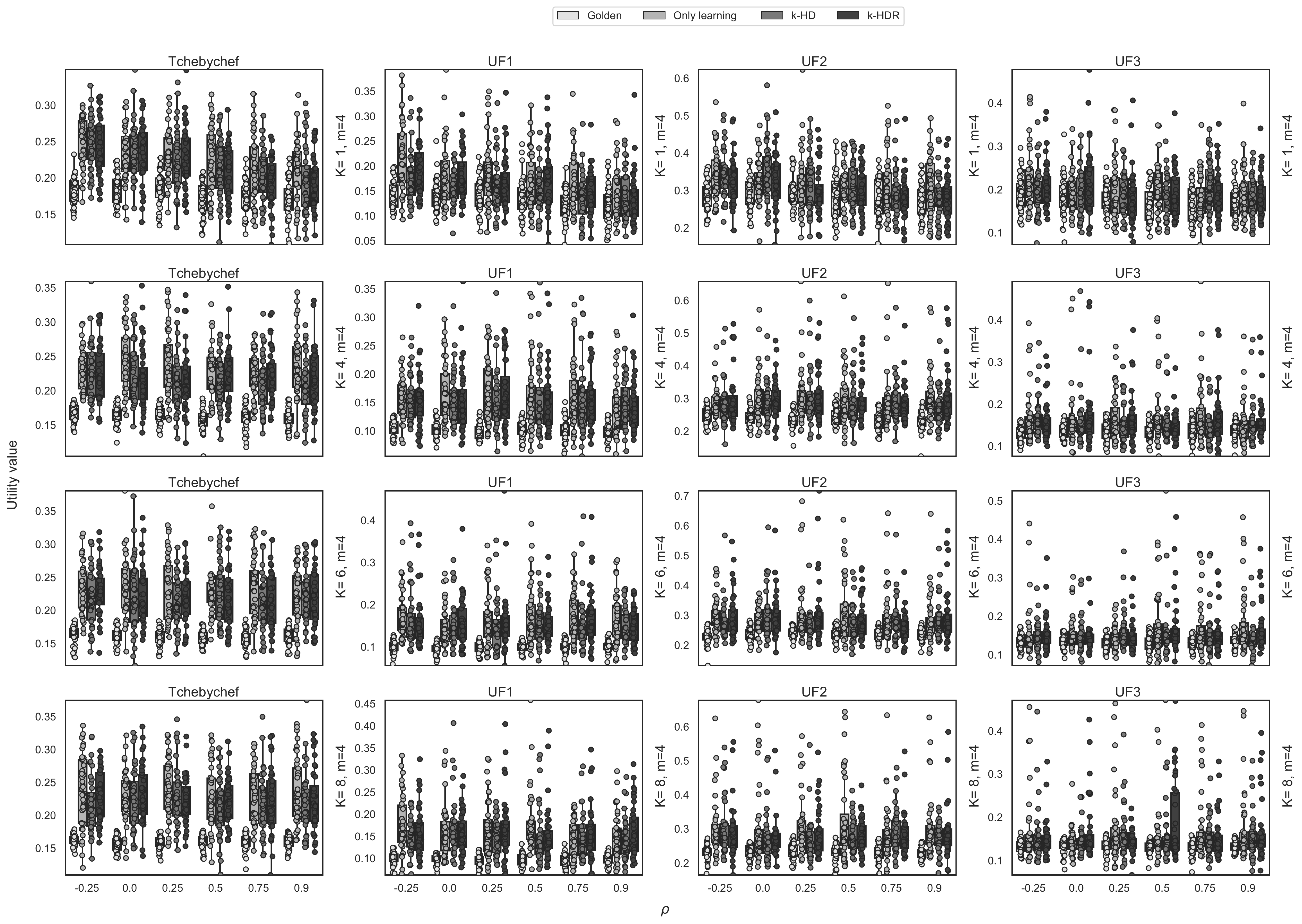}
 \caption{Comparison of the performance of different modes for $\rho$MNK problems with $m=4$. Number of active objectives is fixed.
 }
 \label{fig:rmnkFixed4}
\end{figure*}

\begin{figure*}[!t]
    \includegraphics[width=0.8\textwidth]{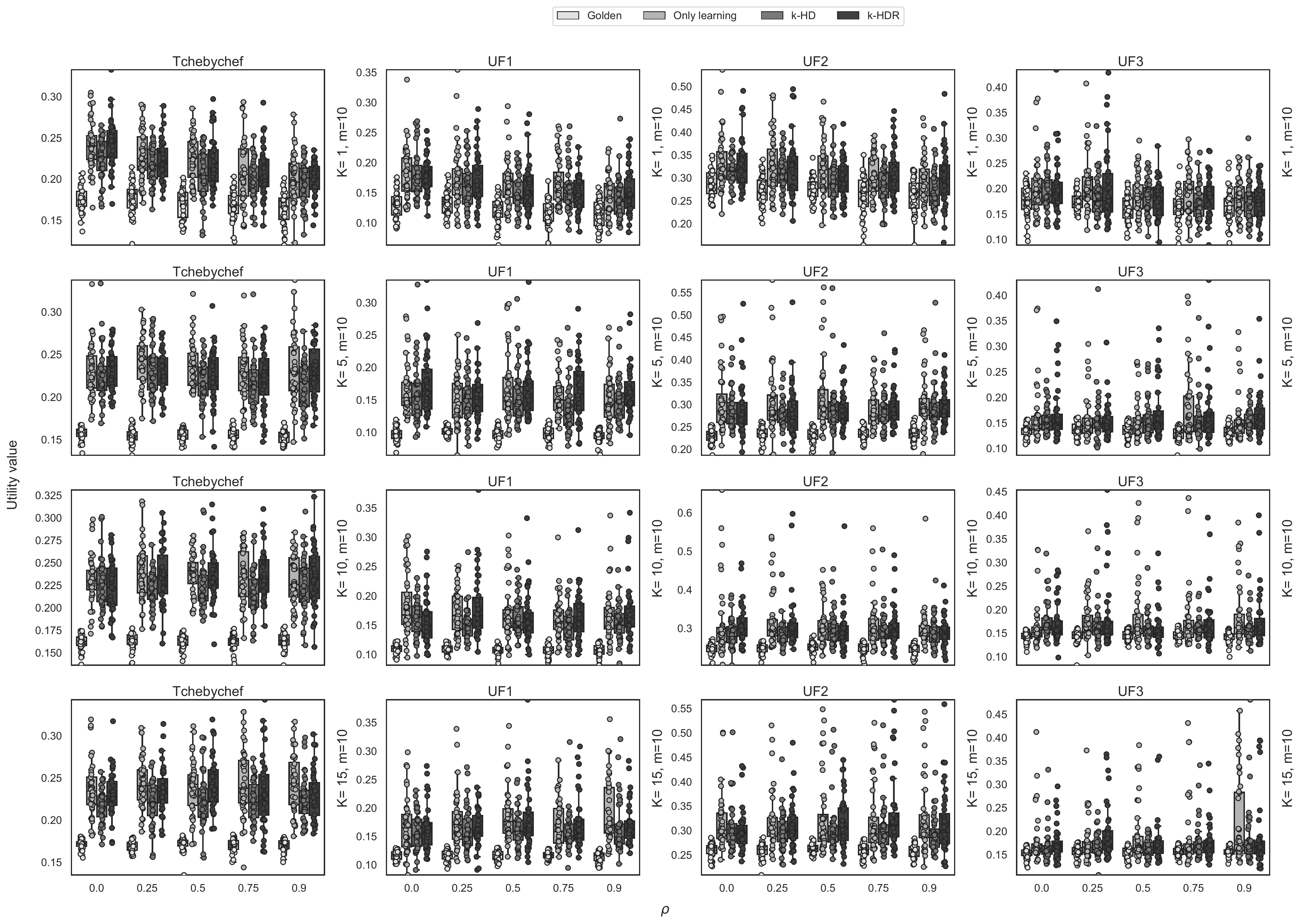}
 \caption{Comparison of the performance of different modes for $\rho$MNK problems with $m=10$. Number of active objectives is fixed.
 }
 \label{fig:rmnkFixed10}
\end{figure*}
\begin{figure*}[!t]
    \includegraphics[width=0.8\textwidth]{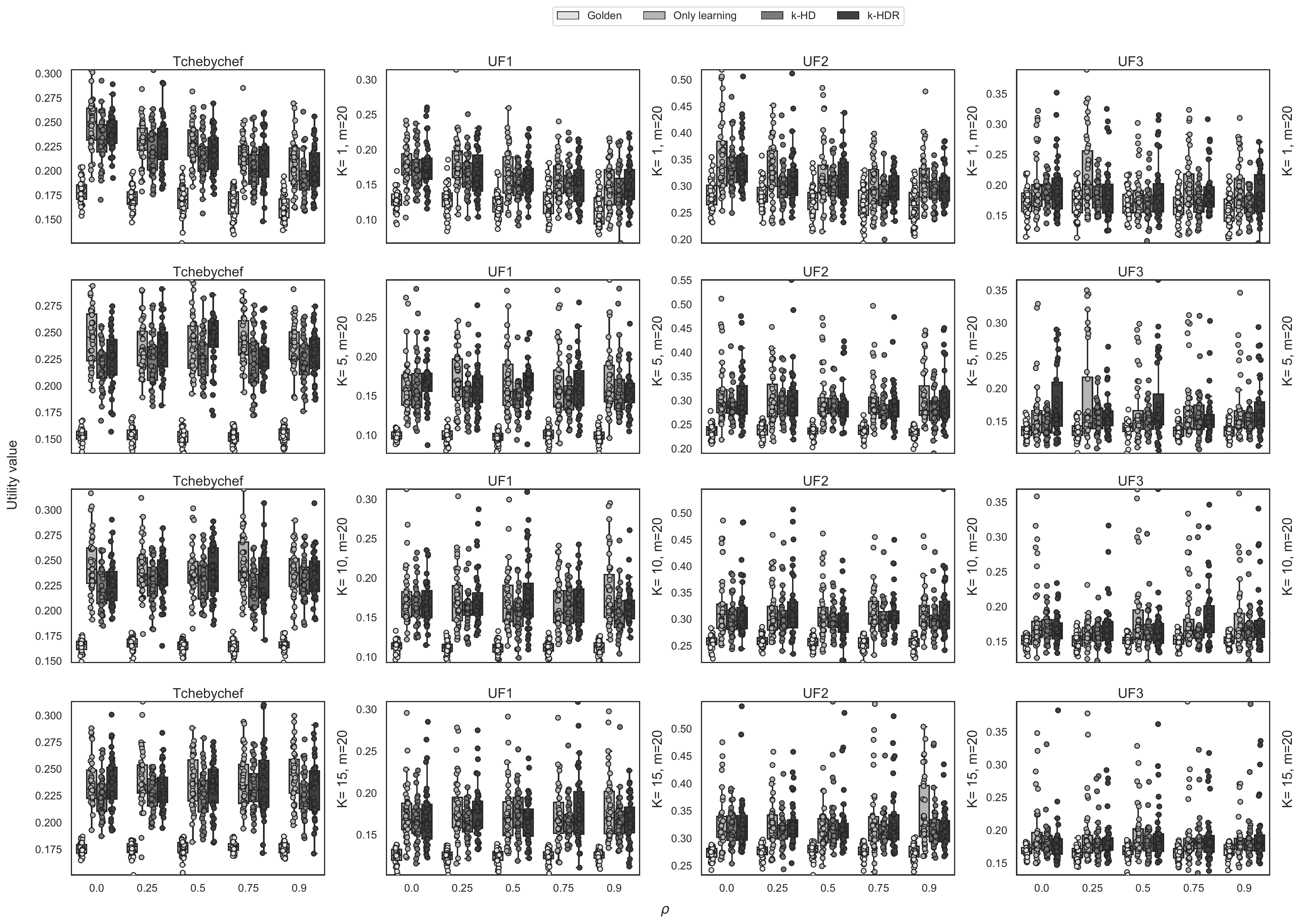}
 \caption{Comparison of the performance of different modes for $\rho$MNK problems with $m=20$. Number of active objectives is fixed.
 }
 \label{fig:rmnkFixed20}
\end{figure*}

\begin{figure*}[!t]
    \includegraphics[width=0.8\textwidth]{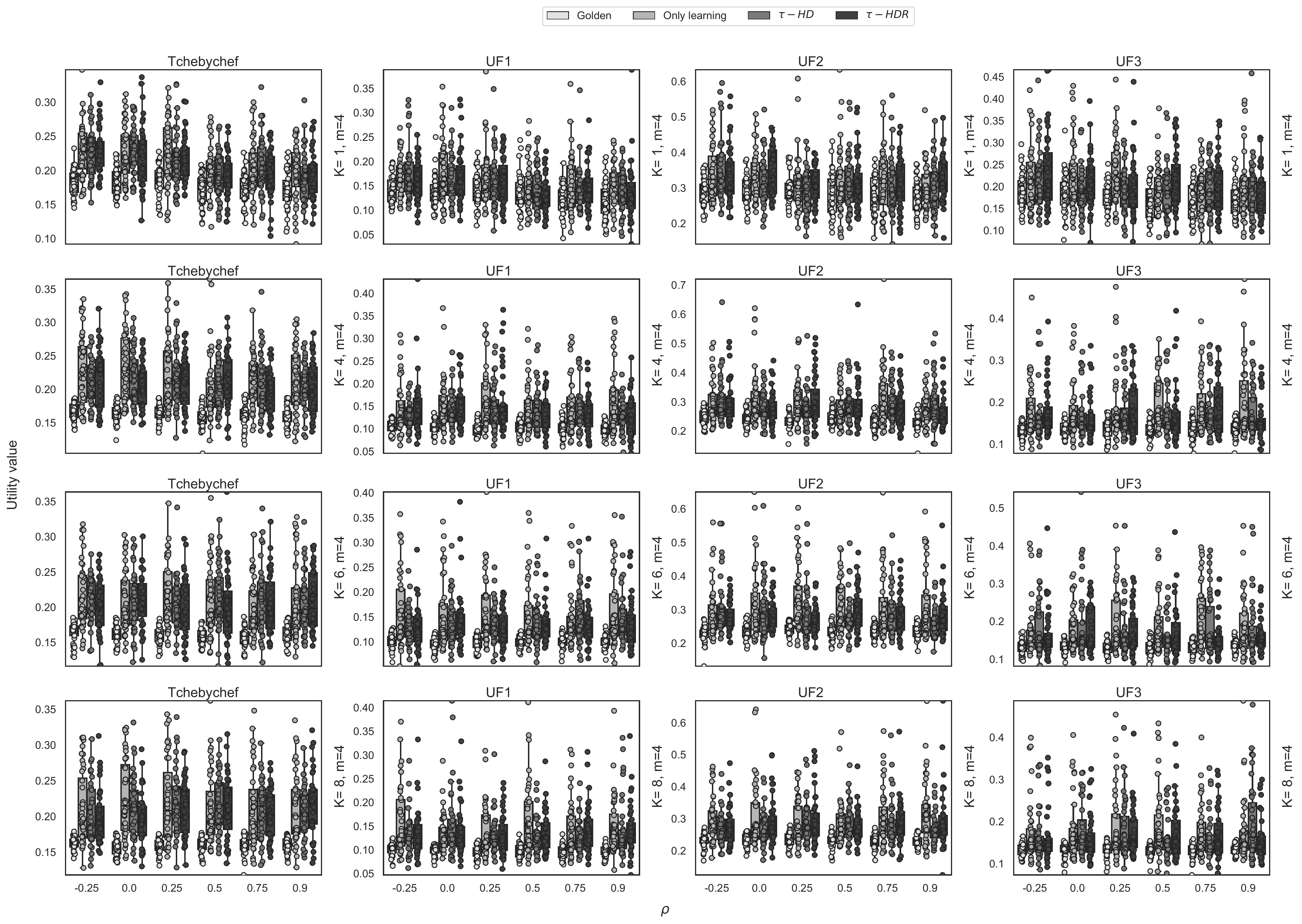}
 \caption{Comparison of the performance of different modes for $\rho$MNK problems with $m=4$. The number of active objectives is not fixed and detection mode is used as a mean of objective reduction technique.
 }
 \label{fig:rmnkreduction4}
\end{figure*}

\begin{figure*}[!t]
    \includegraphics[width=0.8\textwidth]{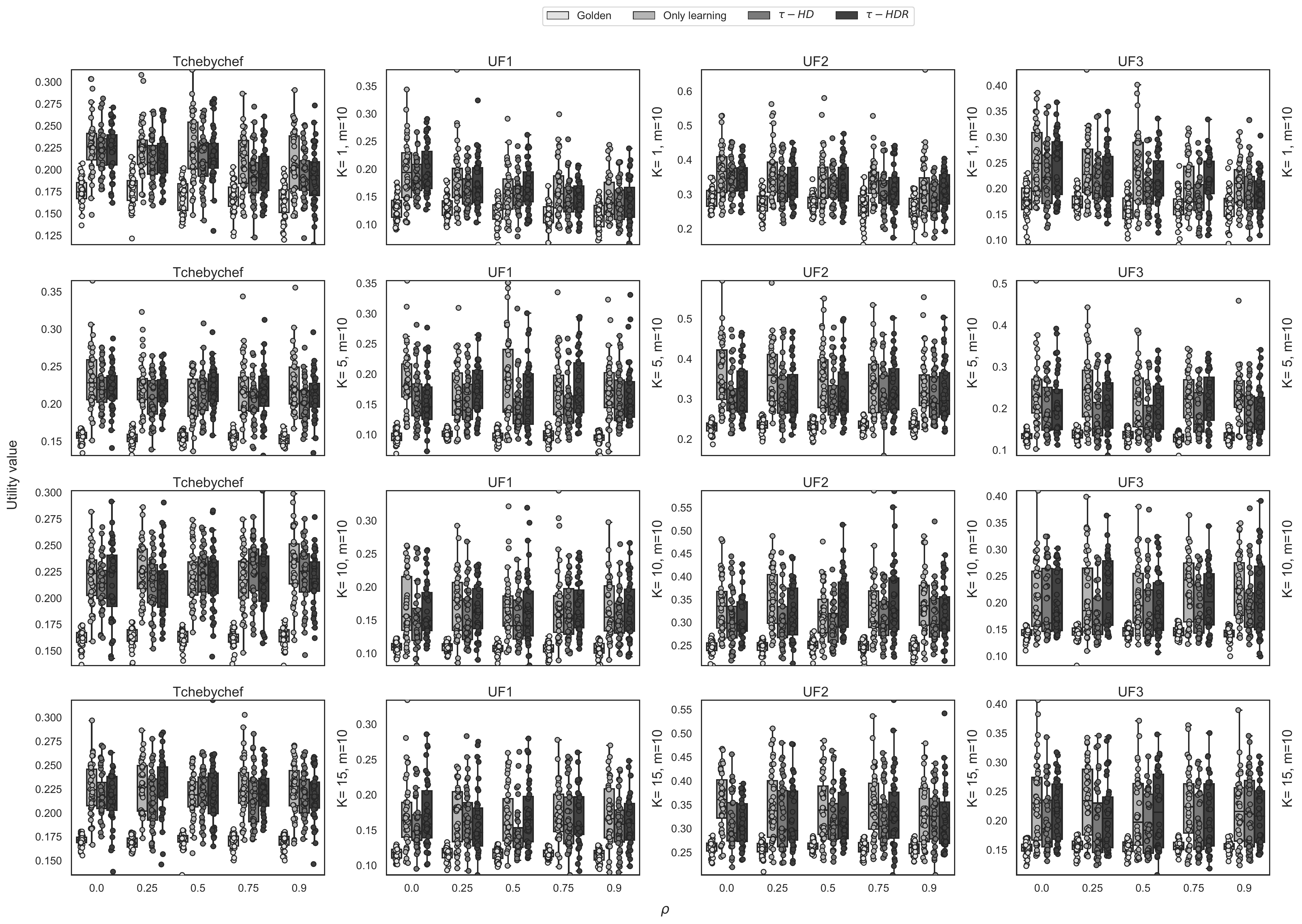}
 \caption{Comparison of the performance of different modes for $\rho$MNK problems with $m=4$. The number of active objectives is not fixed and detection mode is used as a mean of objective reduction technique.
 }
 \label{fig:rmnkreduction10}
\end{figure*}

\begin{figure*}[!t]
    \includegraphics[width=0.8\textwidth]{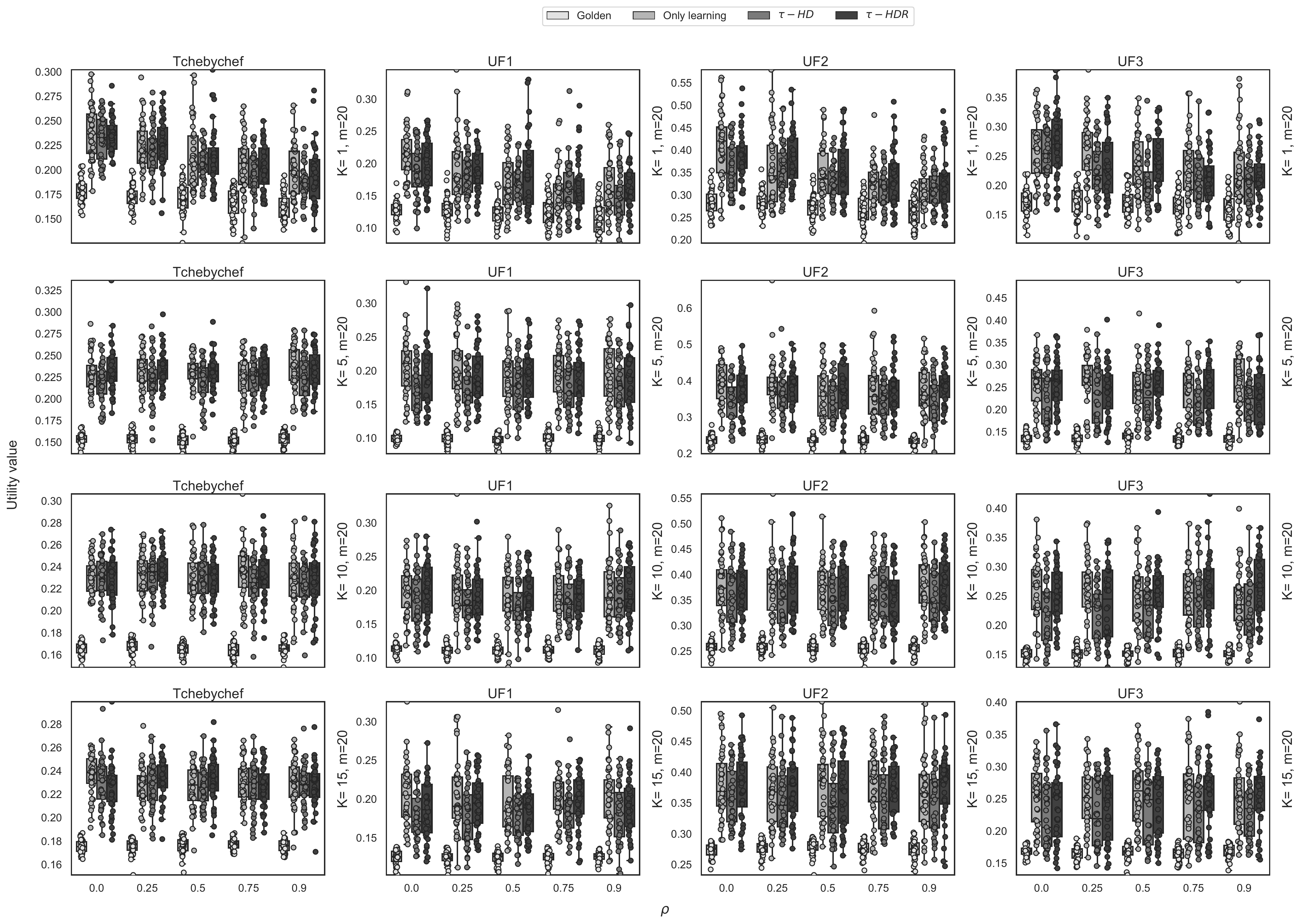}
 \caption{Comparison of the performance of different modes for $\rho$MNK problems with $m=20$. The number of active objectives is not fixed and detection mode is used as a mean of objective reduction technique.
 }
 \label{fig:rmnkreduction20}
\end{figure*}

\begin{figure}[!t]
    \includegraphics[width=0.8\textwidth]{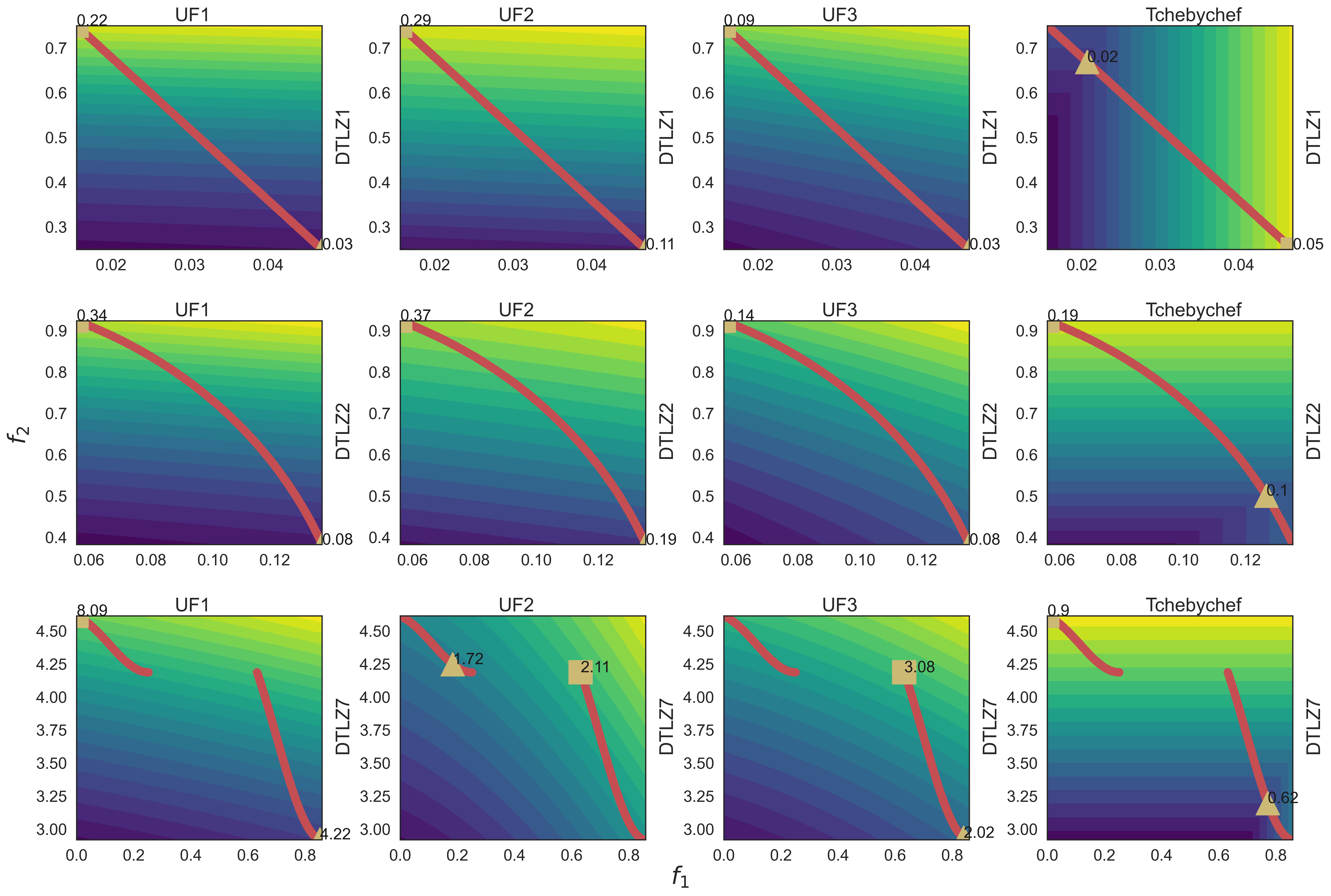}
 of  \caption{PF of the selected problems with first and fourth objectives being active depicted over the contour lines of the different UFs. The position of the worst and best solutions is marked respectively with a  rectangle and a triangle.}
 \label{fig:PF_contour}
 
\end{figure}
\begin{figure*}[!t]
    \includegraphics[width=0.8\textwidth]{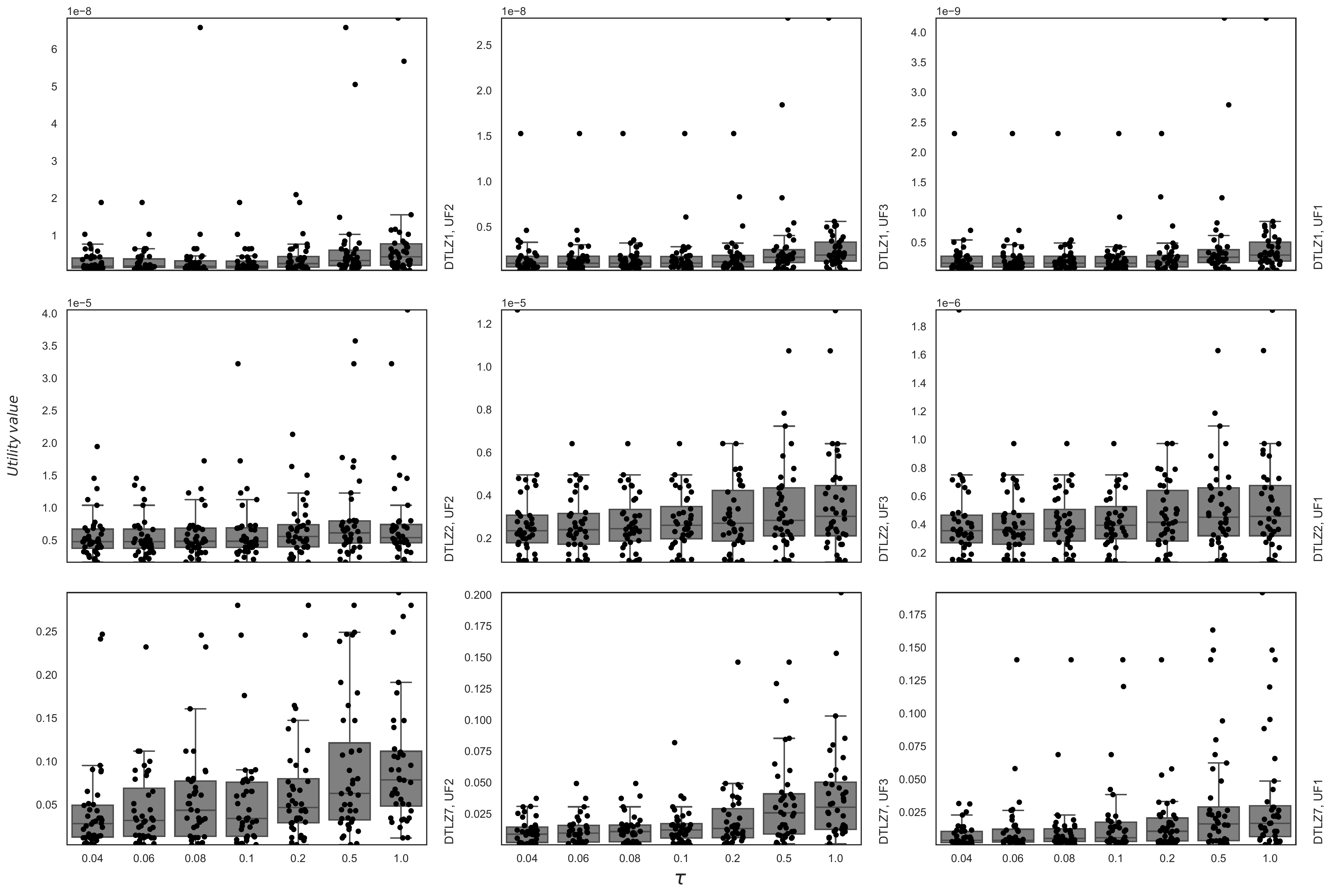}
 \caption{Performance analysis of $\tau$-\emph{HDR} with different values of $\tau$ for DTLZ problems with $m=20$. Number of interaction in all runs is fixed to 6.}
 \label{fig:DTLZ_tau_analysis}
\end{figure*}

\begin{figure*}[!t]
    \includegraphics[width=0.8\textwidth]{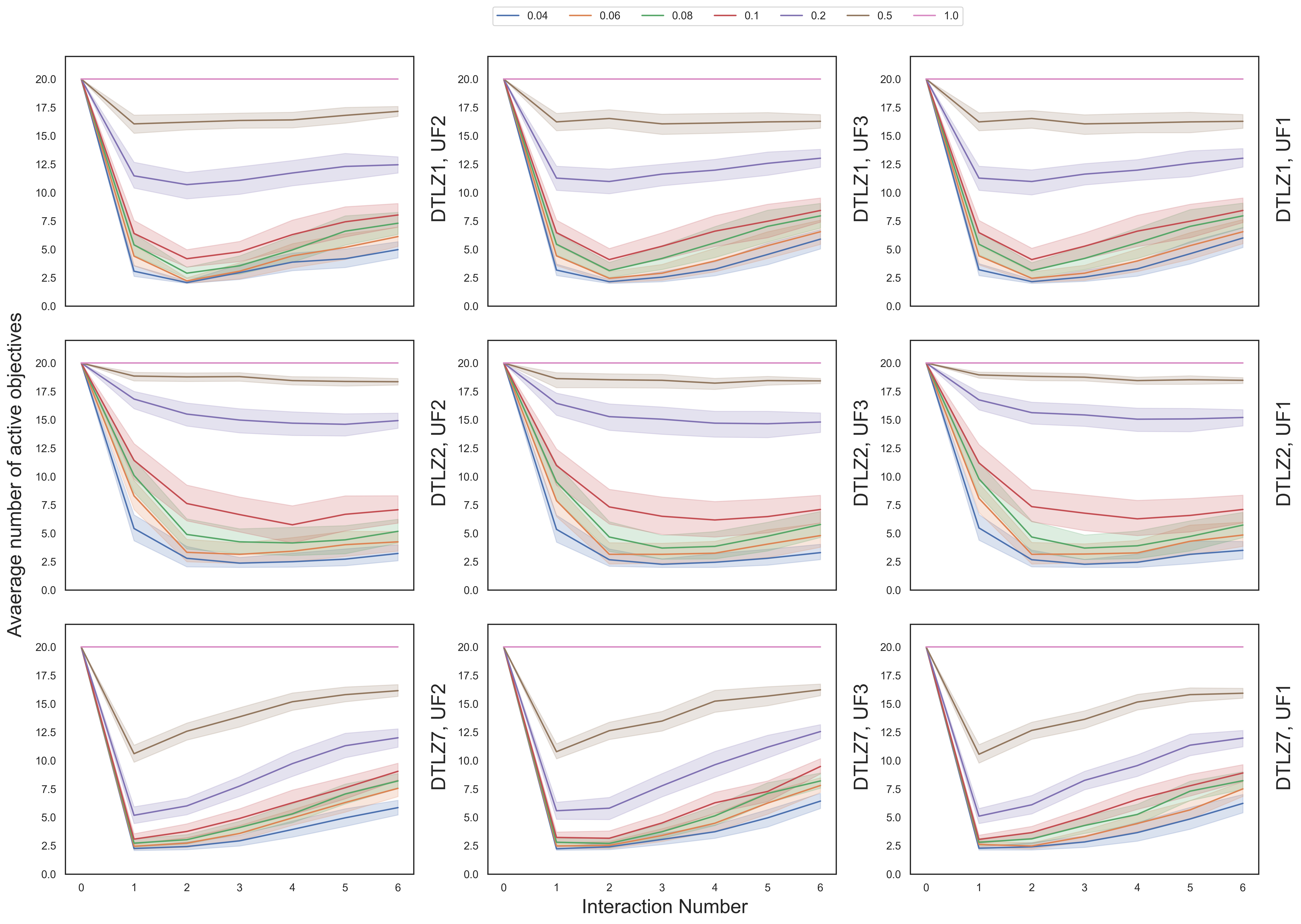}
 \caption{Number of active objectives after each interaction for different values of $\tau$ in $\tau$-\emph{HDR} mode on DTLZ test problems. }
 \label{fig:a_analysis_DTLZ}
\end{figure*}
\begin{figure*}[!t]
    \includegraphics[width=0.8\textwidth]{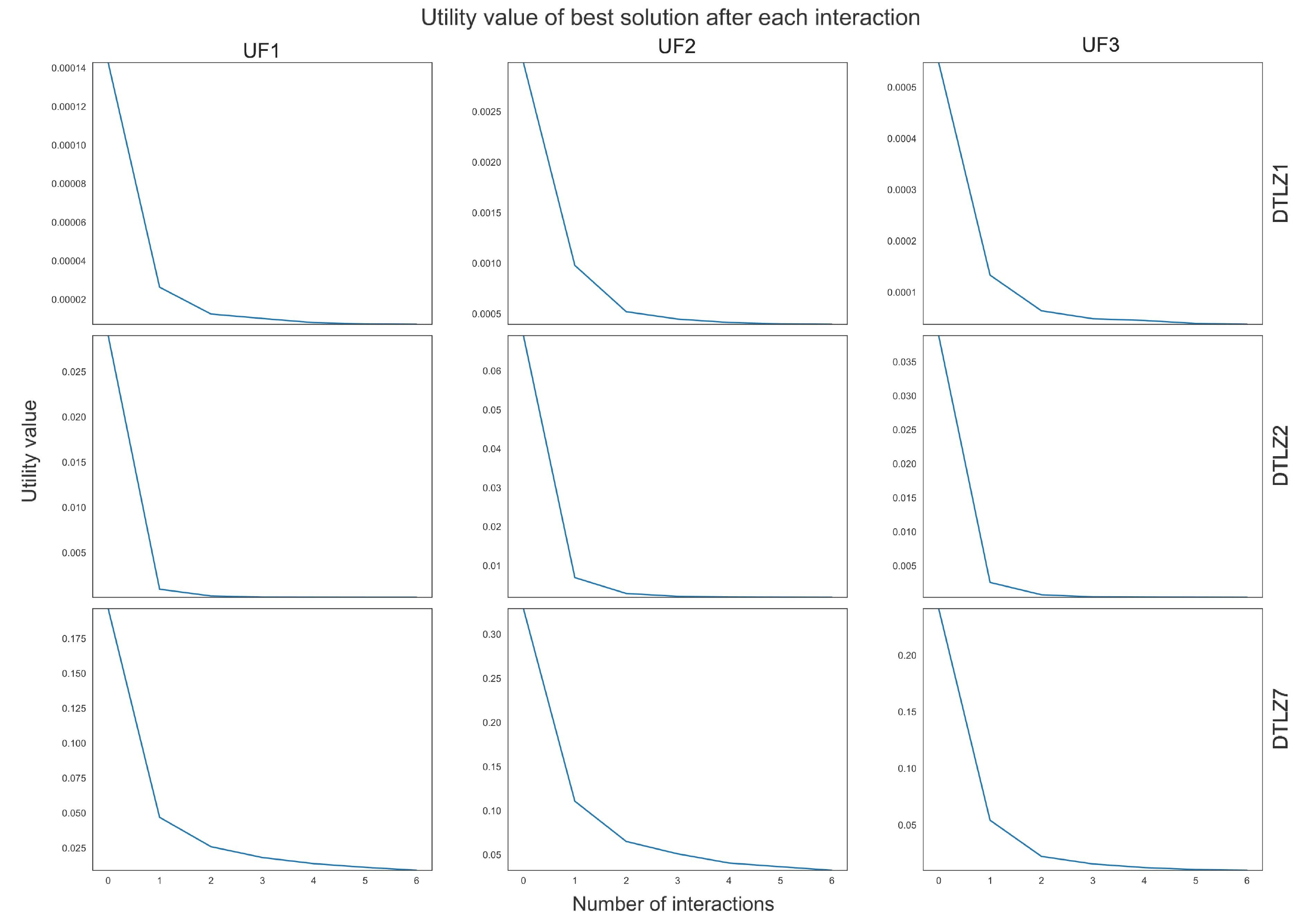}
 \caption{Trend of utility value over different interactions on DTLZ problems when detection mode. }
 \label{fig:trend_dtlz}
\end{figure*}